\documentclass[12pt,reqno]{article}
\usepackage{amssymb}
\usepackage{amsthm}
\usepackage{amsmath}
\newtheorem{theorem}{Theorem}[section]

\newtheorem{definition}[theorem]{Definition}
\newtheorem{example}[theorem]{Example}

\newtheorem{remark}[theorem]{Remark}

\begin{document}
\begin{center} \Large{\bf Invariants of third-order ordinary differential equations $y'''=f(x,y,y',y'')$
via fiber preserving transformations}
\end{center}
\medskip

\begin{center}
Ahmad Y. Al-Dweik$^*$, M. T. Mustafa$^{**}$, H. Azad$^*$ and
F. M. Mahomed$^{***,}$\footnote{Author FM is visiting professorial fellow at UNSW, Sydney for 2014}\\

{$^*$Department of Mathematics \& Statistics, King Fahd University
of Petroleum and Minerals, Dhahran 31261, Saudi Arabia}\\
{$^{**}$Department of Mathematics, Statistics and Physics, Qatar
University, Doha, 2713, State of Qatar}\\
{$^{***}$School of Computational and Applied Mathematics, DST-NRF Centre of Excellence
in Mathematical and Statistical Sciences;
Differential Equations, Continuum Mechanics and Applications,
University of the Witwatersrand, Johannesburg, Wits 2050,  South Africa and\\
$^1$School of Mathematics and Statistics,
University of New South Wales,
Sydney NSW 2052 Australia
}\\

aydweik@kfupm.edu.sa, tahir.mustafa@qu.edu.qa,
hassanaz@kfupm.edu.sa and Fazal.Mahomed@wits.ac.za
\end{center}
\begin{abstract}
Bagderina \cite{Bagderina2008} solved the equivalence problem for
scalar third-order ordinary differential equations (ODEs), quadratic
in the second-order derivative, via point transformations.
However, the question is open for the general class
$y'''=f(x,y,y',y'')$ which is not quadratic in the second-order
derivative. We utilize Lie's infinitesimal method to
study the differential invariants of this general class under
pseudo-group of fiber preserving equivalence transformations
$\bar{x}=\phi(x), \bar{y}=\psi(x,y)$. As a result, all third-order
differential invariants of this group and the invariant
differentiation operators are determined. This leads to simple
necessary explicit conditions for a third-order ODE
to be equivalent to the respective canonical form under the
considered group of transformations. Applications motivated by the literature
are presented.
\end{abstract}
\bigskip
Keywords: Lie's infinitesimal method, differential invariants,
third order ODEs, equivalence problem, fiber preserving
transformations, normal forms, Lie symmetries.
\section{Introduction}
Differential invariants play a significant role in a broad range
of problems arising in  differential geometry, differential
equations, mathematical physics and applications. For instance,
differential invariants have been particularly useful in
dealing with the equivalence problem for geometric structures
\cite{Olver1995}, classification of invariant differential
equations and invariant variational problems
\cite{Kogan2003,Olver1993,Ovsiannikov,Reid2006}, integration of
ordinary differential equations (ODEs) \cite{Olver1993,Ovsiannikov},
equivalence and symmetry properties of solutions \cite{Olver1995},
and in the construction of particular solutions to systems of partial
differential equations (PDEs) \cite{Ovsiannikov,Mansfield2001,Martina2001,Nutku2001}.

Lie \cite{Lie1884} showed that every invariant system of
differential equations \cite{Lie1888} and variational
problem \cite{Lie1897} can be directly expressed in terms of
differential invariants. Along with an illustration \cite{Lie1888}
of how differential invariants could be used to integrate ODEs,
Lie succeeded in completely classifying
all the differential invariants for all possible
finite-dimensional Lie groups of point transformations in the complex plane.

Tress\'e \cite{Tresse1894} and Ovsiannikov  \cite{Ovsiannikov}
gave prominence to Lie's preliminary results on invariant
differentiations and the existence of finite bases of differential
invariants. For the general theory of differential invariants of
Lie groups including algorithms of construction of differential
invariants, the interested reader is referred to
\cite{Olver1995,Ovsiannikov}. Ibragimov \cite{Ibr1997,Ibr1999}
developed a simple method for constructing invariants of families
of linear and nonlinear differential equations admitting infinite
equivalence transformation groups. Lie's infinitesimal method was
applied to solve the equivalence problem for several linear and
nonlinear differential equations
\cite{Ibr2002-1,Ibr2002-2,Ibr2007,Ibr2004,Ibr2005,Johnpillai2001,Torrisi2004,Torrisi2005,Tracina2004}.
Cartan's equivalence method \cite{Olver1995,Gardner1989} is
another systematic approach to solve the equivalence problem for
differential equations. The linearization problem is a particular
case of the equivalence problem.

Linearization criteria for a third-order ODE which are at most
cubic  in the second-order derivative
\begin{equation}
\begin{array}{ll}
y'''=a(x,y,y'){y''}^3+b(x,y,y'){y''}^2+c(x,y,y')y''+d(x,y,y')\label{cub}\\
\end{array}
\end{equation}
have been obtained in \cite{Chern1940,Neut2002} by Cartan's method and then
in \cite{Ibr2005} by the direct approach. Lie \cite{Lie1896} in fact
was the first to note that the third-order ODE
connected via contact transformations to the simplest linear third-order ODE
is at most cubic in the second-order derivative of the form given above in (\ref{cub}).
Recently, Bagderina \cite{Bagderina2009}
presented the basis of differential invariants under the group of contact transformations for the family of
ODEs at most cubic in the second-order derivative (\ref{cub}) by using Lie's method.
She also provided the operators of invariant differentiations.

By using Lie's infinitesimal method, Bagderina \cite{Bagderina2008}
solved the equivalence problem of third-order ODEs which are
at most quadratic  in the second-order derivative
\begin{equation}
\begin{array}{ll}
y'''=a(x,y,y'){y''}^2+b(x,y,y')y''+c(x,y,y')\\
\end{array}
\end{equation}
with respect to the group of point equivalence transformations
\begin{equation}\label{o2}
\begin{array}{ll}
\bar{x}=\phi(x,y), \bar{y}=\psi(x,y).\\
\end{array}
\end{equation}
As an extension, in this paper, we use Lie's infinitesimal method
to study the differential invariants of the third-order ODEs
\begin{equation}\label{w1}
\begin{array}{ll}
y'''=f(x,y,y',y''),\\
\end{array}
\end{equation}
which are not quadratic  in the second-order derivative, under
pseudo-group of fiber preserving equivalence transformations
\begin{equation}\label{o3}
\begin{array}{ll}
\bar{x}=\phi(x), \bar{y}=\psi(x,y).\\
\end{array}
\end{equation}

The structure of the paper is as follows. In the next section, we
give a short description of Lie's infinitesimal method to find the
differential invariants and invariant differentiation operators of
the class of ODEs (\ref{w1}) with respect to the general group of
point equivalence transformations $\bar{x}=\phi(x,y),
\bar{y}=\psi(x,y)$. In Section 3, using the methods described in
Section 2, first, we recover the infinitesimal point equivalence
transformations. Then we find the third-order differential
invariants and invariant differentiation operators of the class of
ODEs (\ref{w1}), which are not quadratic in the second-order
derivative, under two subgroups of the general group of point
equivalence transformations. In Section 4, we provide illustrative
examples of equations not quadratic in the second-order derivative
taken from \cite{Mahomed1996,Zaitsev2002,Dunajski2013}. This is motivated by
studies of this more general class for symmetry properties in
\cite{Mahomed1996}, exact solutions for certain third-order ODEs
in \cite{Zaitsev2002}, linearization and equivalence under contact transformation
for the class (\ref{cub}) in \cite{Ibr2005,Neut2002,Bagderina2009} as well as in
the definition of certain Einstein-Weyl geometry of hyper-CR type which is of recent
interest in physics \cite{Dunajski2013}. Here we consider examples that are
equivalent under the pseudo-group of fiber preserving equivalence
transformations. Another motivation for the examples considered is
that third-order ODEs possessing the Painlev'e property for
polynomial in its lower order derivatives were also investigated,
see e.g. \cite{Mugan1999,Mugan2002}. Finally the conclusion is
presented.

Throughout this paper, we use the notation
$A=[a_1,a_2,...,a_{n}]$ to express any differential operator
$A=\displaystyle{\sum^{n}_{j=1} a_j \frac{\partial}{\partial
b_j}}$. Also, we denote $y', y''$ by $p, q$, respectively.
\section{Lie's infinitesimal method}
In this section, we  briefly describe the Lie method used
to derive differential invariants using point equivalence
transformations.

Consider the $k$th-order system of
PDEs of $n$ independent variables
 $x = (x^{1}, x^{2}, ... , x^{n})$ and $m$ dependent variables $u = (u^{1}, u^{2},...,
 u^{m})$
 \begin{equation}\label{a1}
 \begin{array}{ll}
 E_{\alpha}(x,u,...,u_{(k)})=0, & \alpha=1,...,m~,
\end{array}
\end{equation}
where $u_{(1)},u_{(2)},...,u_{(k)}$ denote the collections of all
first, second, ..., $k$th-order partial derivatives, i.e.,
$u_{i}^{\alpha}=D_{i}(u^{\alpha}),\; u_{ij}^{\alpha}=D_{j}D_{i}(u^{\alpha})$,...,respectively,
in which the total differentiation operator with respect to $x^{i}$
is
\begin{equation}\label{a2}
\begin{array}{ll}
D_{i}=\frac{\partial}{\partial
x^{i}}+u_{i}^{\alpha}\frac{\partial}{\partial
u^\alpha}+u_{ij}^\alpha\frac{\partial}{\partial
u_{j}^\alpha}+...,&i=1,...,n~,
\end{array}
\end{equation}
with the summation convention adopted for repeated indices.
\begin{definition}\rm
\emph{The Lie-B\"acklund operator} is
\begin{equation}\label{a4}
\begin{array}{ll}
X=\xi^{i}\frac{\partial}{\partial
{x^i}}+\eta^\alpha\frac{\partial}{\partial
u^\alpha}&\xi^i,\eta^\alpha\in \emph{A}~,
\end{array}
\end{equation}
where $A$ is the space of  \emph{differential functions}.

The operator (\ref{a4}) is an abbreviated form of the infinite formal
sum
\begin{equation}\label{a5}
\begin{array}{ll}
 X^{(s)}&=\xi^{i}\frac{\partial}{\partial x^{i}}+\eta^{\alpha}\frac{\partial}{\partial u^{\alpha}}+\sum\limits_{s\geq 1 }\zeta_{i_{1}i_{2}...i_{s}}^{\alpha}\frac{\partial}{\partial u_{i_1 i_2...i_s}^\alpha},\\
  &=\xi^{i}D_i+W^\alpha\frac{\partial}{\partial u^{\alpha}}+\sum\limits_{s\geq 1 }D_{i_1}...D_{i_s}(W^\alpha)\frac{\partial}{\partial u_{i_1 i_2...i_s}^\alpha},\\
\end{array}
\end{equation}
where the additional coefficients are determined uniquely by the
prolongation formulae
\begin{equation}\label{a6}
\begin{array}{ll}
\zeta_{i}^\alpha=D_{i}(\eta^{\alpha})-u_{j}^\alpha D_{i}(\xi^{j})=D_{i}(W^\alpha)+\xi^{j}u_{ij}^\alpha,\\
\zeta_{i_1...i_s}^\alpha=D_{i_s}(\zeta_{i_1...i_{s-1}}^\alpha)-u_{ji_1...i_{s-1}}^\alpha D_{i_s}(\xi^{j})=D_{i_1}...D_{i_s}(W^\alpha)+\xi^j u_{ji_1...i_s}^\alpha,& s>1,\\
\end{array}
\end{equation}
in which $W^\alpha$ is the\emph{ Lie characteristic function}
\begin{equation}\label{a7}
W^\alpha=\eta^\alpha-\xi^j u_{j}^{\alpha}.
\end{equation}
\end{definition}
\begin{definition}\rm
\emph{The point equivalence transformation of a class of PDEs
(\ref{a1})} is an invertible transformation of the independent and
dependent variables of the form
\begin{equation}\label{a7}
\begin{array}{ll}
\bar{x}=\phi(x,u), \bar{u}=\psi(x,u),\\
\end{array}
\end{equation}
that maps every equation of the class into an equation of the
same family, viz.
\begin{equation}\label{a8}
\begin{array}{ll}
E_{\alpha}(\bar{x},\bar{u},...,\bar{u}_{(k)})=0, & \alpha=1,...,m.
\end{array}
\end{equation}
\end{definition}
In order to describe Lie's infinitesimal method for deriving
differential invariants using point equivalence transformations,
we use as example the class of equations (\ref{w1}). It is well-known that the
point equivalence transformation
\begin{equation}\label{a9}
\begin{array}{ll}
\bar{x}=\phi(x,y), \bar{y}=\psi(x,y),\\
\end{array}
\end{equation}
 maps (\ref{w1}) into the same family, viz.
\begin{equation}\label{a10}
\bar{y}'''=\bar{f}(\bar{x},\bar{y},\bar{y}',\bar{y}''),
\end{equation}
for arbitrary functions  $\phi(x,y)$ and $\psi(x,y)$, where
$\bar{f}$, in general, can be different from the original function
$f$. The set of all equivalence transformations forms a group
denoted by $\mathcal{E}$.

The standard procedure for Lie's infinitesimal invariance criterion
\cite{Ovsiannikov} is implemented in the next section to recover
the continuous group of point equivalence transformations (\ref
{a9}) for the class of third-order ODEs (\ref{w1}) with the
corresponding infinitesimal point equivalence transformation operator
\begin{equation}\label{a11}
Y=\xi(x,y)D_x+W\partial_y+D_x(W)\partial_{p}+D_x^2(W)\partial_{q}+\mu(x,y,p,q,f)\partial_f,\\
\end{equation}
where  $\xi(x,y)$ and $\eta(x,y)$ are arbitrary functions obtained from
\begin{equation}\label{a12}
\bar{x}= x+\epsilon \xi(x,y)+O(\epsilon^2)=\phi(x,y),\\
\end{equation}
\begin{equation}\label{a13}
\bar{y}= y+\epsilon \eta(x,y)+O(\epsilon^2)=\psi(x,y) ,\\
\end{equation}
and
\begin{equation}\label{a14}
\mu=\dot{D}_x^3(W)+\xi(x,y)\dot{D}_x f,\\
\end{equation} with
$W=\eta-\xi p$ and $ \dot{D}_x = \frac{\partial } {{\partial x}} +
p\frac{\partial }{{\partial y}} + q\frac{\partial } {{\partial p}}
+f\frac{\partial} {{\partial q}}$.
\begin{definition}\rm
\emph{An invariant of a class of third-order ODEs (\ref{w1})} is a
function of the form
 \begin{equation}\label{a15}
\begin{array}{ll}
J=J(x,y,p,q,f),
\end{array}
\end{equation}
which is invariant under the equivalence transformation
(\ref{a9}).
\end{definition}
\begin{definition}\rm
\emph{A differential invariant of order $s$ of a class of
third-order ODEs (\ref{w1})} is a function of the form
 \begin{equation}\label{a16}
\begin{array}{ll}
J=J(x,y,p,q,f,f_{(1)},f_{(2)},...,f_{(s)}),\\
\end{array}
\end{equation}
which is invariant under the equivalence transformation (\ref{a9})
where $f_{(1)},f_{(2)},...,f_{(s)}$ denote the collections of all
first, second,..., $s$th-order partial derivatives.
\end{definition}
\begin{definition}\rm
\emph{An invariant system of order $s$ of a class of third-order
ODEs (\ref{w1})} is the system of the form $E_{\alpha}(x,y,p,q,f,f_{(1)},f_{(2)},...,f_{(s)})=0, ~ \alpha=1,...,m~$ which satisfies the condition\\
 \begin{equation}\label{a22}
\begin{array}{ll}
Y^{(s)}E_{\alpha}(x,y,p,q,f,f_{(1)},f_{(2)},...,f_{(s)})=0~(mod~ E_{\alpha}=0,~\alpha=1,...,m), ~ \alpha=1,...,m.\\
\end{array}
\end{equation}
An invariant system with $\alpha=1$ is called an invariant
equation.
\end{definition}
Now, according to the theory of invariants of infinite
transformation groups \cite{Ovsiannikov}, the invariant criterion
\begin{equation}\label{a17}
Y J(x,y,p,q,f)=0,\\
\end{equation}
should be split by means of the functions  $\xi(x,y)$ and $\eta(x, y)$ and
their derivatives. This gives rise to a homogeneous linear system
of PDEs whose solution gives the required invariants.

It should be noted that since the generator $Y$ contains
arbitrary functions $\xi(x,y)$ and $\eta(x, y)$, the
corresponding identity (\ref{a17}) leads to $m$ linear PDEs for
$J$, where $m$ is the number of the arbitrary functions and their
derivatives that appear in $Y$. We point out that these $m$ PDEs
are not necessarily linearly independent.

In order to determine the differential invariants of order $s$, we
need to calculate the prolongations of the operator $Y$ using
(\ref{a5}) by considering $f$  as a dependent variable and the
variables $x,y,p,q$ as independent variables:
\begin{equation}\label{a18}
\begin{array}{ll}
Y^{(s)}&=Y(x)\tilde{D}_x+Y(y)\tilde{D}_y+Y(p)\tilde{D}_p+Y(q)\tilde{D}_q+\tilde{W}\frac{\partial}{\partial f}+\sum\limits_{s\geq 1 }\tilde{D}_{i_1}...\tilde{D}_{i_s}(\tilde{W})\frac{\partial}{\partial f_{i_1 i_2...i_s}},\\
&i_1, i_2,..., i_s\in \{x,y,p,q\},\\
\end{array}
\end{equation}
where
\begin{equation}\label{a19}
\begin{array}{ll}
\tilde{D}_k= \partial_{k}+f_k \partial_{f}+f_{k i} \partial_{f_{i}}+f_{k i j} \partial_{f_{ij}}+...,& i,j,k\in \{x,y,p,q\}. \\
\end{array}
\end{equation}
in which $\tilde{W}$ is the\emph{ Lie characteristic function}
\begin{equation}\label{a20}
\tilde{W}=\mu-Y(x) f_{x}-Y(y) f_{y}-Y(p) f_{p}-Y(q) f_{q}.
\end{equation}
The differential invariants are determined by the equations
\begin{equation}\label{a21}
Y^{(s)} J(x,y,p,q,f,f_{(1)},f_{(2)},...,f_{(s)})=0.\\
\end{equation}
It should be noted that since the generator $Y^{(s)}$
contains arbitrary functions $\xi(x,y)$ and $\eta(x, y)$,  the
corresponding identity (\ref{a21}) leads to $m$ linear PDEs for
$J$, where $m$ is the number of the arbitrary functions and their
derivatives that appear in $Y^{(s)}$.

For simplicity, from here on, we denote the derivative of
$f(x,y,p,q)$ with respect to the independent variables $x,y,p,q$
 as $f_1,f_2,f_3,f_4$. The same notation will be used for higher-order derivatives.

Now, in order to find all the third order differential invariants
of the third-order ODE (\ref{w1}), one can solve the invariant
criterion (\ref{a21}) with $s=3$. However, for compactness of the
derived differential invariants, one can replace any partial
derivative with respect to $x$ by the total derivative with respect to
$x$. So, we need to solve the following invariant criterion
\begin{small}
\begin{equation}\label{c1}
\begin{array}{ll}
Y^{(3)}J(x,y,p,q,f,f_{{2}},f_{{3}},f_{{4}},f_{{2,2}},f_{{2,3}},f_{{2,4}},f_{{3,3}},f_{{3,4}},f_{{4,4}},f_{{2,2,2}},f_{{2,2,3}},f_{{2,2,4}},f_{{2,3,3}},f_{{2,3,4}},\\
f_{{2,4,4}},f_{{3,3,3}},f_{{3,3,4}},f_{{3,4,4}},f_{{4,4,4}},d_{{1,1}},d_{{1,2}},d_{{1,3}},d_{{1,4}},d_{{1,5}},d_{{1,6}},d_{{1,7}},d_{{1,8}},d_{{1,9}},d_{{1,10}},d_{{2,1}},d_{{2,2}},d_{{2,3}},d_{{2,4}},d_{{3,1}})=0,\\
\end{array}
\end{equation}
\end{small}
by prolonging the infinitesimal operator $Y^{(3)}$ to
the variables $d_{i,j}$ through the infinitesimals
$Y^{(3)}d_{i,j}$, where
\begin{small}
\begin{equation}\label{c2}
\begin{array}{ll}
d_{{1,1}}=\dot{D}_x{f},d_{{1,2}}=\dot{D}_x{f_2},d_{{1,3}}=\dot{D}_x{f_3},d_{{1,4}}=\dot{D}_x{f_4},d_{{1,5}}=\dot{D}_x{f_{2,2}},\\
d_{{1,6}}=\dot{D}_x{f_{2,3}},d_{{1,7}}=\dot{D}_x{f_{2,4}},d_{{1,8}}=\dot{D}_x{f_{3,3}},d_{{1,9}}=\dot{D}_x{f_{3,4}},d_{{1,10}}=\dot{D}_x{f_{4,4}},\\
d_{{2,1}}=\dot{D}^2_x{f},d_{{2,2}}=\dot{D}^2_x{f_2},d_{{2,3}}=\dot{D}^2_x{f_3},d_{{2,4}}=\dot{D}^2_x{f_4},d_{{3,1}}=\dot{D}^3_x{f}.\\
\end{array}
\end{equation}
\end{small}
\begin{definition}\rm
\emph{An invariant differentiation operator of a class of
third-order ODEs (\ref{w1})} is a differential operator
$\mathcal{D}$ which satisfies that if $I$ is a differential
invariant of ODE (\ref{w1}), then $\mathcal{D} I$ is its
differential invariant too.
\end{definition}
As it is shown in \cite{Ovsiannikov}, the number of
independent invariant differentiation operators $\mathcal{D}$
equals  the number of independent variables $x,y,p$ and $q$. The
invariant differentiation operators $\mathcal{D}$ should take the
form
\begin{equation}
\begin{array}{ll}\label{cc}
\mathcal{D}=K\tilde{D}_x+L\tilde{D}_y+M\tilde{D}_p+N\tilde{D}_q,\\
\end{array}
\end{equation}
with the coordinates $K, L, M$ and $N$ satisfying the
non-homogeneous linear system
\begin{equation}
\begin{array}{llll}\label{ccc}
Y^{(3)}{K}=\mathcal{D}(Y(x)),&Y^{(3)}{L}=\mathcal{D}(Y(y)),&Y^{(3)}{M}=\mathcal{D}(Y(p)),&Y^{(3)}{N}=\mathcal{D}(Y(q)),\\
\end{array}
\end{equation}
where $K, L, M$ and $N$ are functions of the following variables
\begin{small}
\begin{equation}\label{cccc}
\begin{array}{ll}
x,y,p,q,f,f_{{2}},f_{{3}},f_{{4}},f_{{2,2}},f_{{2,3}},f_{{2,4}},f_{{3,3}},f_{{3,4}},f_{{4,4}},f_{{2,2,2}},f_{{2,2,3}},f_{{2,2,4}},f_{{2,3,3}},f_{{2,3,4}},f_{{2,4,4}},f_{{3,3,3}},f_{{3,3,4}},\\
f_{{3,4,4}},f_{{4,4,4}},d_{{1,1}},d_{{1,2}},d_{{1,3}},d_{{1,4}},d_{{1,5}},d_{{1,6}},d_{{1,7}},d_{{1,8}},d_{{1,9}},d_{{1,10}},d_{{2,1}},d_{{2,2}},d_{{2,3}},d_{{2,4}},d_{{3,1}}.\\
\end{array}
\end{equation}
\end{small}
In reality, the general solution of the system (\ref{ccc}) gives
both  the differential invariants and the differential
operators. This general solution can be found by prolonging the
infinitesimal operator $Y^{(3)}$ to the variables $K, L,
M$ and $N$ through the infinitesimals $Y^{(3)}{K}, Y^{(3)}{L},
Y^{(3)}{M}$ and $Y^{(3)}{N}$ respectively. Then solving the
invariant criterion
\begin{footnotesize}
\begin{equation}\label{ccccc}
\begin{array}{ll}
Y^{(3)}J(x,y,p,q,f,f_{{2}},f_{{3}},f_{{4}},f_{{2,2}},f_{{2,3}},f_{{2,4}},f_{{3,3}},f_{{3,4}},f_{{4,4}},f_{{2,2,2}},f_{{2,2,3}},f_{{2,2,4}},f_{{2,3,3}},f_{{2,3,4}},f_{{2,4,4}},f_{{3,3,3}},f_{{3,3,4}},\\
f_{{3,4,4}},f_{{4,4,4}},d_{{1,1}},d_{{1,2}},d_{{1,3}},d_{{1,4}},d_{{1,5}},d_{{1,6}},d_{{1,7}},d_{{1,8}},d_{{1,9}},d_{{1,10}},d_{{2,1}},d_{{2,2}},d_{{2,3}},d_{{2,4}},d_{{3,1}},K,L,M,N)=0,\\
\end{array}
\end{equation}
\end{footnotesize}
gives the implicit solution of the variables $K, L, M$ and $N$
with the differential invariants.

In this paper, we are interested in finding the third-order
differential invariants and differential operators of the general
class $y'''=f(x,y,y',y'')$  under a subgroup of point
transformations (\ref{a9}), namely the fiber preserving
transformations $x=\phi(x), y=\psi(x,y)$. So, according to the
theory of invariants of infinite transformation groups
\cite{Ovsiannikov}, the invariant criterion (\ref{ccccc}) should
be split by the functions  $\xi(x)$ and $\eta(x, y)$ and their
derivatives. This gives rise to a homogeneous linear system of
partial differential equations (PDEs):
\begin{equation}
\begin{array}{ll}\label{c3}
X_i J=0, i=1\dots28,&T_i J=0, i=1\dots7,\\
\end{array}
\end{equation}
where $X_i,  i=1\dots28,$ are the differential operators
corresponding to the coefficients of the following derivatives of
$\eta(x, y)$ up to the sixth order in the invariant criterion
\begin{footnotesize}
\begin{equation}
\begin{array}{ll}\label{c4}
\eta_{{}},\eta_{{1}},\eta_{{2}},\eta_{{1,1}},\eta_{{1,2}},\eta_{{2,2}},\eta_{{1,1,1}},\eta_{{1,1,2}},\eta_{{1,2,2}},\eta_{{2,2,2}},\eta_{{1,1,1,1}},\eta_{{1,1,1,2}},\eta_{{1,1,2,2}},\eta_{{1,2,2,2}},\eta_{{2,2,2,2}},\eta_{{1,1,1,1,1}},\eta_{{1,1,1,1,2}},\\
\eta_{{1,1,1,2,2}},\eta_{{1,1,2,2,2}},\eta_{{1,2,2,2,2}},\eta_{{2,2,2,2,2}},\eta_{{1,1,1,1,1,1}},\eta_{{1,1,1,1,1,2}},\eta_{{1,1,1,1,2,2}},\eta_{{1,1,1,2,2,2}},\eta_{{1,1,2,2,2,2}},\eta_{{1,2,2,2,2,2}},\eta_{{2,2,2,2,2,2}}\\
\end{array}
\end{equation}
\end{footnotesize}
and $T_i, i=1\dots7,$ are the differential operators corresponding
to the coefficients of the following derivatives of $\xi(x)$ up to
the sixth order  in the invariant criterion
\begin{footnotesize}
\begin{equation}
\begin{array}{ll}\label{c5}
\xi,\xi_{{1}},\xi_{{1,1}},\xi_{{1,1,1}},\xi_{{1,1,1,1}},\xi_{{1,1,1,1,1}},\xi_{{1,1,1,1,1,1}}.\\
\end{array}
\end{equation}
\end{footnotesize}
The expressions for the differential operators $X_i,  i=1\dots28$
and  $T_i, i=1\dots7$ are therefore too large  and these are given in
the Appendix A, after relabeling the variables
\begin{small}
\begin{equation}\label{c6}
\begin{array}{ll}
x,y,p,q,f,f_{{2}},f_{{3}},f_{{4}},f_{{2,2}},f_{{2,3}},f_{{2,4}},f_{{3,3}},f_{{3,4}},f_{{4,4}},f_{{2,2,2}},f_{{2,2,3}},f_{{2,2,4}},f_{{2,3,3}},f_{{2,3,4}},f_{{2,4,4}},f_{{3,3,3}},f_{{3,3,4}},\\
f_{{3,4,4}},f_{{4,4,4}},d_{{1,1}},d_{{1,2}},d_{{1,3}},d_{{1,4}},d_{{1,5}},d_{{1,6}},d_{{1,7}},d_{{1,8}},d_{{1,9}},d_{{1,10}},d_{{2,1}},d_{{2,2}},d_{{2,3}},d_{{2,4}},d_{{3,1}},K, L, M, N,\\
\end{array}
\end{equation}
\end{small}
by the variables $z_i, i=1\dots43$, respectively.

Functionally independent solutions of system (\ref{c3}) provide
all independent differential invariants of $y'''=f(x,y,y',y'')$ up
to the third order under the fiber preserving transformation, as
well as an implicit solution of the variables $K, L, M$ and $N$
which yield the differential operators via (\ref{cc})  as
explained in the next section.

The existence of the solutions for system (\ref{c3}) is proved by
showing that the differential operators $X_i,  i=1\dots28$ and
$T_i, i=1\dots7$ form a Lie algebra $\mathcal{L}_{35}$
\cite[p.422, Theorem 14.1]{Olver1995}. The nonzero commutators for
the Lie algebra $\mathcal{L}_{35}$ are given in the Appendix B,
after relabeling the differential operators $X_i,  i=1\dots28$ and
$T_i, i=1\dots7$ by the operators $e_i, i=1\dots35$, respectively.

Using Appendix B, it can be seen that the  Lie algebra
$\mathcal{L}_{35}$ is solvable and has the chain of Lie
subalgebras $0=\mathcal{G}_{0} \lhd\mathcal{G}_{1}
\lhd\mathcal{G}_{2} \lhd\mathcal{G}_{3}
\lhd\mathcal{G}_{4}\lhd\mathcal{G}_{5}\lhd\mathcal{G}_{6}\lhd\mathcal{G}_{7}\lhd\mathcal{G}_{8}=\mathcal{L}_{35}$
with each an ideal in the next where
\begin{small}
\begin{equation}\label{cad}
\begin{array}{ll}
\mathcal{G}_{1}= \{e_{22}, e_{23}, e_{24}, e_{25}, e_{26}, e_{27}, e_{28}\},\\
\mathcal{G}_{2} = \mathcal{G}_{1}  \cup \{e_{16}, e_{17}, e_{18}, e_{19}, e_{20}, e_{21}\},\\
\mathcal{G}_{3} = \mathcal{G}_{2}  \cup \{e_{11}, e_{12}, e_{13}, e_{14}, e_{15}\},\\
\mathcal{G}_{4}= \mathcal{G}_{3}   \cup \{e_{7}, e_{8}, e_9, e_{10}\},\\
\mathcal{G}_{5}= \mathcal{G}_{4}   \cup \{e_4, e_5, e_6\},\\
\mathcal{G}_{6}= \mathcal{G}_{5}   \cup \{e_{33},e_{34},e_{35}\},\\
\mathcal{G}_{7}= \mathcal{G}_{6}   \cup \{e_{31},e_{32}\},\\
\mathcal{G}_{8}= \mathcal{G}_{7}   \cup \{e_1,e_2,e_3,e_{29},e_{30}\}.\\
\end{array}
\end{equation}
\end{small}
In the next section, we solve the system (\ref{c3}) by using the
chain (\ref{cad}). In more detail, as $\mathcal{G}_{1}$ is
abelian, one can find its joint invariants by finding the
invariants of its generators in any order. Since $\mathcal{G}_{1}$
is an ideal in $\mathcal{G}_{2}$, the algebra $\mathcal{G}_{2}$
operates on the joint invariants of $\mathcal{G}_{1}$. Moreover,
as $\mathcal{G}_{1}$ is abelian, then the induced action of
$\mathcal{G}_{2}$ on the joint invariants of $\mathcal{G}_{1}$ is
also abelian. Continuing in this manner, we obtain the joint
invariants of the full algebra.

The reader is further referred to \cite{Ibr1997} and \cite[Section
10]{Ibr1999} for examples illustrating the infinitesimal method.
\section{Third-order differential invariants and invariant equations for $y''' = f(x,y,y',y'')$}

\subsection{The infinitesimal point equivalence transformations}
In order to find continuous group of equivalence transformations
of the class (\ref{w1}) we consider the arbitrary function $f$
that appears in our equation as a dependent variable and the
variables $x,y,y'=p,y''=q$ as independent variables and apply the
Lie infinitesimal invariance criterion \cite{Ovsiannikov}, that is
we look for the infinitesimal $\xi, \eta$ and $\mu$ of the
equivalence operator $Y$:
\begin{equation}
Y=\xi(x,y)\partial_x+\eta(x,y)\partial_y+\mu(x,y,p,q,f)\partial_f,\\
\end{equation}
such that its prolongation leaves the equation (\ref{w1})
invariant.

The prolongation of operator $Y$ can be given using (\ref{a5}) as
\begin{equation}
Y=\xi(x,y)D_x+W\partial_y+D_x(W)\partial_{p}+D_x^2(W)\partial_{q}+D_x^3(W)\partial_{y'''}+\mu(x,y,p,q,f)\partial_f,\\
\end{equation}
where $$ D_x  = \frac{\partial } {{\partial x}} + p\frac{\partial
}{{\partial y}} + q\frac{\partial } {{\partial p}}
+y'''\frac{\partial} {{\partial q}}+y^{(4)}\frac{\partial}
{{\partial y'''}}+...$$ is the operator of total derivative and $W
= \eta(x,y) - \xi(x,y) p$ is the characteristic of infinitesimal
operator $X=\xi(x,y)\partial_x+\eta(x,y)\partial_y$.

So, the Lie infinitesimal invariance criterion gives
$\mu=\dot{D}_x^3(W)+\xi(x,y)\dot{D}_x f$ for arbitrary functions
$\xi(x,y)$ and $\eta(x,y)$ where $ \dot{D}_x = \frac{\partial }
{{\partial x}} + p\frac{\partial }{{\partial y}} + q\frac{\partial
} {{\partial p}} +f\frac{\partial} {{\partial q}}$.

Thus, equation (\ref{w1}) admits an infinite continuous group of
equivalence transformations generated by the Lie algebra
$\mathcal{L}_\mathcal{E}$ spanned by the following infinitesimal
operators
\begin{equation}\label{e3}
U=\xi(x,y)\frac{\partial } {{\partial x}}-p D_x(\xi)\partial_{p}-(2q D_x(\xi)+p D_x^2(\xi))\partial_{q}-(3 f D_x(\xi)+3 q D_x^2(\xi)+p \dot{D}_x^3(\xi))\partial_f,\\
\end{equation}
\begin{equation}\label{e4}
V=\eta(x,y)\partial_y+D_x(\eta)\partial_{p}+D_x^2(\eta)\partial_{q}+\dot{D}_x^3(\eta)\partial_f,\\
\end{equation}
The infinitesimal point equivalence transformations
(\ref{e3})-(\ref{e4}) can be written in the finite form as in
(\ref{a12})-(\ref{a13}), respectively, where $\phi$ and $\psi$ are
arbitrary functions of the indicated variables.
\subsection{Third-order differential invariants and invariant equations under the transformation $\bar{x}=x,~\bar{y}=\psi(x,y)$}
In this section, we derive all the third order differential
invariants of the general class $y'''=f(x,y,y',y'')$  under a
subgroup of point transformations (\ref{a9}), namely the
transformations $\bar{x}=x,~\bar{y}=\psi(x,y)$. Moreover, the
invariant differentiation operators \cite{Ovsiannikov} are
constructed in order to get some higher-order differential
invariants from the lower-order ones. Precisely, we obtain the
following theorem.
\begin{theorem}
Let $y'''=f(x,y,y',y'')$ be the  class of third-order ODE with
$f_{4,4,4}\ne0$. All the third order differential invariants,
under the point transformations $\bar{x}=x,~ \bar{y}=\psi(x,y)$, are
functions of the following seventeen differential invariants
\begin{equation}\label{d1}
\begin{array}{llllll}
\alpha_{{1}}=x,&\alpha_{{2}}={\frac {\lambda_{{5}}}{\lambda_{{4}}}},&\alpha_{{3}}={\frac {\lambda_{{6}}}{{\lambda_{{4}}}^{2}}},&\alpha_{{4}}={\frac {\lambda_{{7}}}{{\lambda_{{4}}}^{2}}},&\alpha_{{5}}={\frac {\lambda_{{8}}}{{\lambda_{{4}}}^{2}}} ,&\alpha_{{6}}={\frac {\lambda_{{9}}}{{\lambda_{{4}}}^{2}}},\\
\alpha_{{7}}={\frac {\lambda_{{10}}}{{\lambda_{{4}}}^{2}}},&\alpha_{{8}}={\frac {\lambda_{{11}}}{{\lambda_{{4}}}^{2}}},&\alpha_{{9}}=\lambda_{{12}},&\alpha_{{10}}=\lambda_{{13}},&\alpha_{{11}}={\frac {\lambda_{{14}}}{\lambda_{{4}}}},&\alpha_{{12}}={\frac {\lambda_{{15}}}{\lambda_{{4}}}},\\
\alpha_{{13}}={\frac{\lambda_{{16}}}{\lambda_{{4}}}},&\alpha_{{14}}={\frac{\lambda_{{17}}}{\lambda_{{4}}}},&\alpha_{{15}}={\frac {\lambda_{{18}}}{\lambda_{{4}}}},&\alpha_{{16}}=\lambda_{{19}},&\alpha_{{17}}=\lambda_{{20}},&\\
\end{array}
\end{equation}
where $\{\lambda_{{i}}\}^{20}_{i=4}$ are the following relative
invariants
\begin{footnotesize}
\begin{equation*}
\begin{array}{lllll}
\lambda_{{4}}=&f_{{4,4}},\\
\lambda_{{5}}=&\frac{1}{3}\left(2\,f_{{3,4}}f_{{4}}-2f_{{3}}f_{{4,4}}-6f_{{2,4}}+3f_{{3,3}}\right),\\
\lambda_{{6}}=&f_{{4,4,4}},\\
\lambda_{{7}}=&\frac{1}{3}\left(2f_{{4}}f_{{4,4,4}}+3f_{{3,4,4}}\right),\\
\lambda_{{8}}=&\frac{1}{9}\left(4f_{{4,4,4}}f^2_{{4}}+12f_{{4}}f_{{3,4,4}}+4f_{{4}}{f^2_{{4,4}}}+9f_{{3,3,4}}+6f_{{3,4}}f_{{4,4}}\right),\\
\lambda_{{9}}=&\frac{1}{9}\left(2f_{{4,4,4}}{f^2_{{4}}}+2f_{{4}}{f^2_{{4,4}}}+3f_{{4}}f_{{3,4,4}}+3f_{{4,4,4}}f_{{3}}+9f_{{2,4,4}}+3f_{{3,4}}f_{{4,4}}\right),\\
\lambda_{{10}}=&\frac{1}{18}(-4{f^2_{{4}}}f_{{3,4,4}}+4f_{{4}}f_{{4,4,4}}f_{{3}}+12f_{{4}}f_{{2,4,4}}-12f_{{3,3,4}}f_{{4}}-4f_{{4}}f_{{3,4}}f_{{4,4}}-6{f^2_{{3,4}}}-9f_{{3,3,3}}\\
&+4{f^2_{{4,4}}}f_{{3}}+6f_{{3}}f_{{3,4,4}}+12f_{{2,4}}f_{{4,4}}+18f_{{2,3,4}}) ,\\
\lambda_{{11}}=&\frac{1}{27}(-4f_{{4}}f_{{3}}{f^2_{{4,4}}}-6f_{{4,4}}f_{{3,4}}f_{{3}}-18f_{{4,4}}f_{{2,3}}-12f_{{4,4}}f_{{2,4}}f_{{4}}+4f_{{4,4}}f_{{3,4}}{f^2_{{4}}}-6f_{{4,4,4}}{f^2_{{3}}}-36f_{{3}}f_{{2,4,4}}\\
&-4f_{{3}}f_{{4,4,4}}{f^2_{{4}}}+9f_{{3}}f_{{3,3,4}}+12{f^2_{{3,4}}}f_{{4}}+9f_{{3,4}}f_{{3,3}}+4{f^3_{{4}}}f_{{3,4,4}}-12{f^2_{{4}}}f_{{2,4,4}}+12f_{{3,3,4}}{f^2_{{4}}}\\
&+9f_{{4}}f_{{3,3,3}}-54f_{{2,2,4}}+27f_{{2,3,3}}),\\
\lambda_{{12}}=&\frac{1}{3}(-{f^2_{{4}}}-3f_{{3}}+3\dot{D}_x{f_4}),\\
\end{array}
\end{equation*}
\end{footnotesize}
\begin{footnotesize}
\begin{equation}\label{d2}
\begin{array}{lllll}
\lambda_{{13}}=&\frac{1}{9}(-2{f^3_{{4}}}-9f_{{4}}f_{{3}}+6f_{{4}}\dot{D}_x{f_4}-27f_{{2}}+9\dot{D}_x{f_3}) ,\\
\lambda_{{14}}=&\frac{1}{3}(3\dot{D}_x{f_{4,4}}+f_{{4,4}}f_{{4}}),\\
\lambda_{{15}}=&\frac{1}{9}(2f_{{4,4}}{f^2_{{4}}}+3f_{{3}}f_{{4,4}}+6f_{{4}}\dot{D}_x{f_{4,4}}-9f_{{2,4}}+9\dot{D}_x{f_{3,4}}),\\
\lambda_{{16}}=&\frac{1}{9}(4{f^2_{{4}}}\dot{D}_x{f_{4,4}}-2{f^2_{{4}}}f_{{3,4}}-3f_{{4}}f_{{3,3}}+12f_{{4}}\dot{D}_x{f_{3,4}}+4f_{{4}}f_{{4,4}}\dot{D}_x{f_4}-12f_{{2,4}}f_{{4}}+6f_{{3,4}}\dot{D}_x{f_4}-18f_{{2,3}}+9\dot{D}_x{f_{3,3}}),\\
\lambda_{{17}}=&\frac{1}{27}(2f_{{4,4}}{f^3_{{4}}}-6{f^2_{{4}}}f_{{3,4}}+6{f^2_{{4}}}\dot{D}_x{f_{4,4}}-9f_{{4}}f_{{3,3}}+12f_{{4}}f_{{3}}f_{{4,4}}+9f_{{4}}\dot{D}_x{f_{3,4}}-9f_{{2,4}}f_{{4}}+9f_{{3}}\dot{D}_x{f_{4,4}}-27f_{{2,3}}\\
&+27f_{{2}}f_{{4,4}}+27\dot{D}_x{f_{2,4}}),\\
\lambda_{{18}}=&\frac{1}{27}(-9{f^2_{{3}}}f_{{4,4}}+9f_{{3}}\dot{D}_x{f_{3,4}}+6f_{{3}}f_{{4}}\dot{D}_x{f_{4,4}}+6f_{{3}}f_{{4,4}}\dot{D}_x{f_4}+6f_{{3}}f_{{3,4}}f_{{4}}-54f_{{3}}f_{{2,4}}+9f_{{3}}f_{{3,3}}+4f_{{4,4}}{f^2_{{4}}}\dot{D}_x{f_4}\\
&+6f_{{3,4}}f_{{4}}\dot{D}_x{f_4}+18f_{{2,4}}\dot{D}_x{f_4}+27f_{{2}}f_{{3,4}}-2f_{{3,4}}{f^3_{{4}}}-30{f^2_{{4}}}f_{{2,4}}+18f_{{4}}f_{{4,4}}f_{{2}}+4{f^3_{{4}}}\dot{D}_x{f_{4,4}}+12{f^2_{{4}}}\dot{D}_x{f_{3,4}}\\
&-3{f^2_{{4}}}f_{{3,3}}+9f_{{4}}\dot{D}_x{f_{3,3}}+18f_{{4}}\dot{D}_x{f_{2,4}}-45f_{{4}}f_{{2,3}}-81f_{{2,2}}+27\dot{D}_x{f_{2,3}}),\\
\lambda_{{19}}=&\frac{1}{9}(-2{f^3_{{4}}}-9f_{{4}}f_{{3}}-27f_{{2}}+9\dot{D}^2_x{f_4}),\\
\lambda_{{20}}=&\frac{1}{27}(-2{f^4_{{4}}}-12{f^2_{{4}}}f_{{3}}-6{f^2_{{4}}}\dot{D}_x{f_4}+18f_{{4}}\dot{D}^2_x{f_4}-27f_{{4}}\dot{D}_x{f_3}-18{f^2_{{3}}}+9f_{{3}}\dot{D}_x{f_4}-81\dot{D}_x{f_2}+27\dot{D}^2_x{f_3}) ,\\
\end{array}
\end{equation}
\end{footnotesize}
Moreover, the invariant differential operators are
\begin{footnotesize}
\begin{equation}\label{d3}
\begin{array}{llll}
\mathcal{D}_1=\frac{1}{f_{4,4}}\tilde{D}_q,\\
\mathcal{D}_2=\frac{1}{f_{4,4}}(3\tilde{D}_p+2f_4\tilde{D}_q),\\
\mathcal{D}_3=\frac{1}{f_{4,4}}(9\tilde{D}_y+3f_4\tilde{D}_p+(2{f_4}^2+3f_3)\tilde{D}_q),\\
\mathcal{D}_4=\tilde{D}_x+p\tilde{D}_y+q\tilde{D}_p+f\tilde{D}_q.\\
\end{array}
\end{equation}
\end{footnotesize}
\end{theorem}
\proof Functionally independent solutions of the subsystem
\begin{equation}\label{s1}
\begin{array}{ll}
X_i J=0, i=1\dots28,\\
\end{array}
\end{equation}
of (\ref{c3}) provide all independent differential invariants of
$y'''=f(x,y,y',y'')$ up to the third order under the
transformations $\bar{x}=x,~\bar{y}=\psi(x,y)$, as well as an
implicit solution of the variables $K, L, M$ and $N$ which provide
the differential operators via (\ref{cc}).

 The solution of system (\ref{s1}) is found in two steps using Maple through the chain (\ref{cad}). First we
consider the following subsystem of equations (\ref{s1})
\begin{equation}\label{s2}
X_i J=0, i=4\dots28.
\end{equation}
In 43-dimensional space of the variables $z_i, i=1\dots43$, the rank
of the system (\ref{s2}) is 19, so it has 24 functionally
independent solutions which are given as:
\begin{footnotesize}
\begin{equation}
\begin{array}{lllll}
\lambda_{{1}}=&z_{{1}},\\
\lambda_{{2}}=&z_{{2}},\\
\lambda_{{3}}=&z_{{3}}, \\
\lambda_{{4}}=&z_{{14}},\\
\lambda_{{5}}=&\frac{2}{3}\,z_{{13}}z_{{8}}-\frac{2}{3}\,z_{{7}}z_{{14}}-2\,z_{{11}}+z_{{12}},\\
\lambda_{{6}}=&z_{{24}},\\
\lambda_{{7}}=&\frac{2}{3}\,z_{{8}}z_{{24}}+z_{{23}},\\
\lambda_{{8}}=&\frac{4}{9}\,z_{{24}}{z_{{8}}}^{2}+\frac{4}{9}\,z_{{8}}{z_{{14}}}^{2}+\frac{4}{3}\,z_{{23}}z_{{8}}+z_{{22}}+\frac{2}{3}\,z_{{13}}z_{{14}},\\
\lambda_{{9}}=&\frac{2}{9}\,z_{{24}}{z_{{8}}}^{2}+\frac{1}{3}\,z_{{23}}z_{{8}}+\frac{2}{9}\,z_{{8}}{z_{{14}}}^{2}+\frac{1}{3}\,z_{{24}}z_{{7}}+z_{{20}}+\frac{1}{3}\,z_{{13}}z_{{14}}\\
\lambda_{{10}}=&-\frac{2}{9}\,{z_{{8}}}^{2}z_{{23}}-\frac{2}{3}\,z_{{22}}z_{{8}}+\frac{2}{9}\,z_{{8}}z_{{24}}z_{{7}}+\frac{2}{3}\,z_{{8}}z_{{20}}-\frac{2}{9}\,z_{{8}}z_{{13}}z_{{14}}-\frac{1}{3}\,{z_{{13}}}^{2}-\frac{1}{2}\,z_{{21}}+\frac{1}{3}\,z_{{7}}z_{{23}}+\frac{2}{9}\,{z_{{14}}}^{2}z_{{7}}\\
&+\frac{2}{3}\,z_{{14}}z_{{11}}+z_{{19}} ,\\
\lambda_{{11}}=&{\frac {4}{27}}\,{z_{{8}}}^{3}z_{{23}}-{\frac{4}{27}} \,{z_{{8}}}^{2}z_{{24}}z_{{7}}+{\frac{4}{27}}\,{z_{{8}}}^{2}z_{{13}}z_{{14}}+\frac{4}{9}\,z_{{22}}{z_{{8}}}^{2}-\frac{4}{9}\,{z_{{8}}}^{2}z_{{20}}-{\frac{4}{27}}\,{z_{{14}}}^{2}z_{{8}}z_{{7}}-\frac{4}{9}\,z_{{8}}z_{{14}}z_{{11}}\\
&+\frac{4}{9}\,{z_{{13}}}^{2}z_{{8}}+\frac{1}{3}\,z_{{21}}z_{{8}}-\frac{2}{9}\,z_{{13}}z_{{7}}z_{{14}}+\frac{1}{3}\,z_{{12}}z_{{13}}-\frac{2}{9}\,z_{{24}}{z_{{7}}}^{2}+\frac{1}{3}\,z_{{7}}z_{{22}}-\frac{4}{3}\,z_{{7}}z_{{20}}-2\,z_{{17}}-\frac{2}{3}\,z_{{14}}z_{{10}}+z_{{18}},\\
\lambda_{{12}}=&-\frac{1}{3}\,{z_{{8}}}^{2}-z_{{7}}+z_{{28}},\\
\lambda_{{13}}=&-\frac{2}{9}\,{z_{{8}}}^{3}-z_{{7}}z_{{8}}+\frac{2}{3}\,z_{{8}}z_{{28}}-3\,z_{{6}}+z_{{27}} ,\\
\lambda_{{14}}=&\frac{1}{3}\,z_{{14}}z_{{8}}+z_{{34}},\\
\lambda_{{15}}=&\frac{2}{9}\,z_{{14}}{z_{{8}}}^{2}+\frac{2}{3}\,z_{{8}}z_{{34}}+\frac{1}{3}\,z_{{7}}z_{{14}}-z_{{11}}+z_{{33}},\\
\lambda_{{16}}=&\frac{4}{9}\,{z_{{8}}}^{2}z_{{34}}-\frac{2}{9}\,{z_{{8}}}^{2}z_{{13}}+\frac{4}{9}\,z_{{8}}z_{{14}}z_{{28}}+\frac{4}{3}\,z_{{8}}z_{{33}}-\frac{4}{3}\,z_{{11}}z_{{8}}-\frac{1}{3}\,z_{{8}}z_{{12}}+\frac{2}{3}\,z_{{28}}z_{{13}}-2\,z_{{10}}+z_{{32}},\\
\lambda_{{17}}=&{\frac{2}{27}}\,z_{{14}}{z_{{8}}}^{3}+\frac{4}{9}\,z_{{14}}z_{{7}}z_{{8}}+z_{{6}}z_{{14}}-\frac{2}{9}\,{z_{{8}}}^{2}z_{{13}}+\frac{2}{9}\,{z_{{8}}}^{2}z_{{34}}+\frac{1}{3}\,z_{{8}}z_{{33}}-\frac{1}{3}\,z_{{11}}z_{{8}}-\frac{1}{3}\,z_{{8}}z_{{12}}+\frac{1}{3}\,z_{{7}}z_{{34}}-z_{{10}}+z_{{31}},\\
\lambda_{{18}}=&{\frac {4}{27}}\,{z_{{8}}}^{3}z_{{34}}-{\frac{2}{27}} \,{z_{{8}}}^{3}z_{{13}}-{\frac{10}{9}}\,{z_{{8}}}^{2}z_{{11}}+{\frac{4}{27}}\,{z_{{8}}}^{2}z_{{14}}z_{{28}}-\frac{1}{9}\,{z_{{8}}}^{2}z_{{12}}+\frac{4}{9}\,{z_{{8}}}^{2}z_{{33}}\\
&+\frac{2}{9}\,z_{{8}}z_{{13}}z_{{7}}+\frac{2}{3}\,z_{{8}}z_{{6}}z_{{14}}+\frac{1}{3}\,z_{{8}}z_{{32}}+\frac{2}{3}\,z_{{8}}z_{{31}}-\frac{5}{3}\,z_{{8}}z_{{10}}+\frac{2}{9}\,z_{{8}}z_{{7}}z_{{34}}+\frac{2}{9}\,z_{{8}}z_{{28}}z_{{13}}\\
&+\frac{2}{9}\,z_{{14}}z_{{28}}z_{{7}}-\frac{1}{3}\,z_{{14}}{z_{{7}}}^{2}+\frac{2}{3}\,z_{{28}}z_{{11}}+z_{{30}}-3\,z_{{9}}-2\,z_{{7}}z_{{11}}+z_{{6}}z_{{13}}+\frac{1}{3}\,z_{{7}}z_{{33}}+\frac{1}{3}\,z_{{7}}z_{{12}},\\
\lambda_{{19}}=&-\frac{2}{9}\,{z_{{8}}}^{3}-z_{{7}}z_{{8}}-3\,z_{{6}}+z_{{38}},\\
\lambda_{{20}}=&-{\frac{2}{27}}\,{z_{{8}}}^{4}-\frac{4}{9}\,{z_{{8}}}^{2}z_{{7}}-\frac{2}{9}\,{z_{{8}}}^{2}z_{{28}}-z_{{8}}z_{{27}}+\frac{2}{3}\,z_{{8}}z_{{38}}-\frac{2}{3}\,{z_{{7}}}^{2}+\frac{1}{3}\,z_{{28}}z_{{7}}-3\,z_{{26}}+z_{{37}} ,\\
\end{array}
\end{equation}
\end{footnotesize}
and
\begin{footnotesize}
\begin{equation}\label{c8}
\begin{array}{lllll}
\lambda_{{21}}=&z_{{40}},\\
\lambda_{{22}}=&z_{{41}},\\
\lambda_{{23}}=&\frac{1}{3}\,z_{{8}}z_{{40}}z_{{3}}-\frac{1}{3}\,z_{{41}}z_{{8}}-z_{{40}}z_{{4}}+z_{{42}} ,\\
\lambda_{{24}}=&\frac{2}{3}\,z_{{40}}z_{{8}}z_{{4}}+\frac{1}{3}\,z_{{40}}z_{{7}}z_{{3}}-z_{{40}}z_{{5}}-\frac{1}{3}\,z_{{7}}z_{{41}}-\frac{2}{3}\,z_{{8}}z_{{42}}+z_{{43}},\\
\end{array}
\end{equation}
\end{footnotesize}
In variables $\lambda_i, i= i=1\dots24$, the remaining non-zero
operators in system (\ref{s1}) take the form
\begin{small}
\begin{equation}\label{c9}
\begin{array}{ll}
X_1=&[0,1,0,0,0,0,0,0,0,0,0,0,0,0,0,0,0,0,0,0,0,0,0,0],\\
X_2=&[0,0,1,0,0,0,0,0,0,0,0,0,0,0,0,0,0,0,0,0,0,\lambda_{{21}},0,0],\\
X_3=&[0,0,\lambda_{{3}},-\lambda_{{4}},-\lambda_{{5}},-2\,\lambda_{{6}},-2\,\lambda_{{7}},-2\,\lambda_{{8}},-2\,\lambda_{{9}},-2\,\lambda_{{10}},\\
&-2\,\lambda_{{11}},0,0,-\lambda_{{14}},-\lambda_{{15}},-\lambda_{{16}},-\lambda_{{17}},-\lambda_{{18}},0,0,0,\lambda_{{22}},\lambda_{{23}},\lambda_{{24}}],\\
\end{array}
\end{equation}
\end{small}
Finally, we consider the following subsystem of equations
(\ref{c3})
\begin{equation}\label{t13}
X_i J=0, i=1\dots3.
\end{equation}
In 24-dimensional space of variables $\lambda_i, i=1\dots24$, the
rank of the system (\ref{t13}) is 3, so it has 21 functionally
independent solutions which are given as:
\begin{equation}\label{t14}
\begin{array}{llllll}
\alpha_{{1}}=x,&\alpha_{{2}}={\frac {\lambda_{{5}}}{\lambda_{{4}}}},&\alpha_{{3}}={\frac {\lambda_{{6}}}{{\lambda_{{4}}}^{2}}},&\alpha_{{4}}={\frac {\lambda_{{7}}}{{\lambda_{{4}}}^{2}}},&\alpha_{{5}}={\frac {\lambda_{{8}}}{{\lambda_{{4}}}^{2}}} ,&\alpha_{{6}}={\frac {\lambda_{{9}}}{{\lambda_{{4}}}^{2}}},\\
\alpha_{{7}}={\frac {\lambda_{{10}}}{{\lambda_{{4}}}^{2}}},&\alpha_{{8}}={\frac {\lambda_{{11}}}{{\lambda_{{4}}}^{2}}},&\alpha_{{9}}=\lambda_{{12}},&\alpha_{{10}}=\lambda_{{13}},&\alpha_{{11}}={\frac {\lambda_{{14}}}{\lambda_{{4}}}},&\alpha_{{12}}={\frac {\lambda_{{15}}}{\lambda_{{4}}}},\\
\alpha_{{13}}={\frac{\lambda_{{16}}}{\lambda_{{4}}}},&\alpha_{{14}}={\frac{\lambda_{{17}}}{\lambda_{{4}}}},&\alpha_{{15}}={\frac {\lambda_{{18}}}{\lambda_{{4}}}},&\alpha_{{16}}=\lambda_{{19}},&\alpha_{{17}}=\lambda_{{20}},&\\
\end{array}
\end{equation}
and
\begin{equation}\label{t15}
\begin{array}{llll}
\alpha_{{18}}=\lambda_{{21}},&\alpha_{{19}}=- \lambda_{{ 4}}\left(\lambda_{{3}}\lambda_{{21}}-\lambda_{{22}} \right) ,&\alpha_{{20}}=\lambda_{{4}}\lambda_{{23}},&\alpha_{{21}}=\lambda_{{4}} \lambda_{{24}}.\\
\end{array}
\end{equation}
Here $\alpha_{{18}}, \alpha_{{19}}, \alpha_{{20}}$ and
$\alpha_{{21}}$ are the only invariants depending on the variables
$K, L, M$ and $N$. Then the general solution of (\ref{ccc}) , for
$\xi=0,\eta=\eta(x,y)$, can be given implicitly by back
substitution as
\begin{equation}\label{t16}
\begin{array}{llll}
K_{{}}={\it F_1},\\
f_{{4,4}}\left( p K_{{}}-L_{{}} \right)={\it F_2},\\
f_{{4,4}}\left(f_4(pK-L)-3(qK-M)\right)={\it F_3} ,\\
f_{{4,4}}\left(f_3(pK-L)+2f_4(qK-M)-3(fK-N)\right)=\,{\it F_4} .\\
\end{array}
\end{equation}
where $F_1, F_2, F_3$ and $F_4$ are the arbitrary functions of
$\alpha_{{i}}, i=1...17.$

Finally, solving system (\ref{t16}) gives the variables $K, L, M$
and $N$ in terms of four arbitrary functions $F_1, F_2, F_3$ and
$F_4$ which provide four independent invariant differentiation
operators $\mathcal{D}_1,\mathcal{D}_2,\mathcal{D}_3$ and
$\mathcal{D}_4$ via (\ref{cc}).
\endproof
\subsection{Third-order differential invariants and invariant equations under the fiber
preserving transformation $\bar{x}=\phi(x), \bar{y}=\psi(x,y)$}
This section is devoted to the derivation of all third order differential
invariants of the general class $y'''=f(x,y,y',y'')$  under a
subgroup of point transformations (\ref{a9}), namely the fiber
preserving transformations $\bar{x}=\phi(x),~\bar{y}=\psi(x,y)$.
The invariant differentiation operators are also constructed.
Precisely, we obtain the following theorem.
\begin{theorem}
Let $y'''=f(x,y,y',y'')$ be the  class of third-order ODE with
$f_{4,4,4}\ne0$. All the third order differential invariants,
under pseudo-group of fiber preserving transformations
$\bar{x}=\phi(x), \bar{y}=\psi(x,y)$, are functions of the
following eleven differential invariants
\begin{equation}\label{d4}
\begin{array}{llllll}
\beta_{{1}}={\frac {\gamma_{{6}}\gamma_{{5}}}{{\gamma_{{4}}}^{4}}},&\beta_{{2}}={\frac{\gamma_{{7}}\gamma_{{5}}}{{\gamma_{{4}}}^{4}}},&\beta_{{3}}={\frac{\gamma_{{8}}{\gamma_{{5}}}^{2}}{{\gamma_{{4}}}^{6} }} ,&\beta_{{4}}={\frac{\gamma_{{9}}{\gamma_{{5}}}^{3}}{{\gamma_{{4}}}^{8} }} ,&\beta_{{5}}={\frac{{\gamma_{{5}}}^{2}\gamma_{{10}}}{{\gamma_{{4}}}^{4 }}} ,&\beta_{{6}}={\frac{\gamma_{{5}}\gamma_{{11}}}{{\gamma_{{4}}}^{3}}},\\
\beta_{{7}}={\frac{\gamma_{{12}}{\gamma_{{5}}}^{2}}{{\gamma_{{4}}}^{5}}},&\beta_{{8}}={\frac{\gamma_{{13}}{\gamma_{{5}}}^{3}}{{\gamma_{{4}}}^{7}}},&\beta_{{9}}={\frac{\gamma_{{14}}{\gamma_{{5}}}^{3}}{{\gamma_{{4}}}^{7}}},&\beta_{{10}}={\frac{\gamma_{{15}}{\gamma_{{5}}}^{4}}{{\gamma_{{4}}}^{9}}},&\beta_{{11}}={\frac{\gamma_{{16}}{\gamma_{{5}}}^{3}}{{\gamma_{{4}}}^{6}}},\\
\end{array}
\end{equation}
where $\{\gamma_{{i}}\}^{16}_{i=4}$ are relative invariants given
by (\ref{c11}).

Moreover, the invariant differential operators are
\begin{footnotesize}
\begin{equation}\label{d5}
\begin{array}{llll}
\mathcal{D}_1=\frac{f_{4,4}}{f_{4,4,4}}\tilde{D}_q,\\
\mathcal{D}_2=\frac{1}{f_{4,4}f_{4,4,4}}\left(f_{4,4,4}\tilde{D}_p-f_{3,4,4}\tilde{D}_q\right),\\
\mathcal{D}_3=\frac{1}{f^4_{4,4}f_{4,4,4}}\left(6f_{4,4}f^2_{4,4,4}\tilde{D}_y-2f_{4,4}f_{4,4,4}(f_{{4}}f_{{4,4,4}}+3f_{{3,4,4}})\tilde{D}_p+(3f_{{4,4}}f^2_{{3,4,4}}+f^2_{{4, 4,4}}\left( 3f_{{3,3}}+2f_{{3,4}}f_{{4}}-6f_{{2,4}} \right))\tilde{D}_q\right),\\
\mathcal{D}_4=\frac{f_{4,4,4}}{f^2_{4,4}}\left(\tilde{D}_x+p\tilde{D}_y+q\tilde{D}_p+f\tilde{D}_q\right).\\
\end{array}
\end{equation}
\end{footnotesize}
\end{theorem}
\proof Functionally independent solutions of the system (\ref{c3})
provide all independent differential invariants of
$y'''=f(x,y,y',y'')$ up to the third order under the fiber
preserving transformations $\bar{x}=\phi(x), \bar{y}=\psi(x,y)$,
as well as an implicit solution of the variables $K, L, M$ and $N$
which provide the differential operators via (\ref{cc}).

In variables $\lambda_i, i= i=1\dots24$, given in (\ref{d2}), the
remaining non-zero operators in system (\ref{c3}) take the form
\begin{small}
\begin{equation}\label{c9}
\begin{array}{ll}
X_1=&[0,1,0,0,0,0,0,0,0,0,0,0,0,0,0,0,0,0,0,0,0,0,0,0],\\
X_2=&[0,0,1,0,0,0,0,0,0,0,0,0,0,0,0,0,0,0,0,0,0,\lambda_{{21}},0,0],\\
X_3=&[0,0,\lambda_{{3}},-\lambda_{{4}},-\lambda_{{5}},-2\,\lambda_{{6}},-2\,\lambda_{{7}},-2\,\lambda_{{8}},-2\,\lambda_{{9}},-2\,\lambda_{{10}},\\
&-2\,\lambda_{{11}},0,0,-\lambda_{{14}},-\lambda_{{15}},-\lambda_{{16}},-\lambda_{{17}},-\lambda_{{18}},0,0,0,\lambda_{{22}},\lambda_{{23}},\lambda_{{24}}],\\
T_1=&[1,0,0,0,0,0,0,0,0,0,0,0,0,0,0,0,0,0,0,0,0,0,0,0],\\
T_2=&[0,0,-\lambda_{{3}},\lambda_{{4}},-\lambda_{{5}},3\,\lambda_{{6}},2\,\lambda_{{7}},\lambda_{{8}},\lambda_{{9}},0,-\lambda_{{11}},-2\,\lambda_{{12}},-3\,\lambda_{{13}},0,\\
&-\lambda_{{15}},-2\,\lambda_{{16}},-2\,\lambda_{{17}},-3\,\lambda_{{18}},-3\,\lambda_{{19}},-4\,\lambda_{{20}},\lambda_{{21}},0,-\lambda_{{23}},-2\,\lambda_{{24}}],\\
T_3=&[0,0,0,0,0,0,-\lambda_{{6}},-\frac{2}{3}\,{\lambda_{{4}}}^{2}-2\,\lambda_{{7}},-\lambda_{{7}}-\frac{1}{3}\,{\lambda_{{4}}}^{2},\frac{1}{2}\,\lambda_{{8}}-\lambda_{{9}},2\,\lambda_{{10}}+\frac{1}{3}\,\lambda_{{4}}\lambda_{{5}},0,-\lambda_{{12}},0,\\
&-\lambda_{{14}},-\frac{2}{3}\,\lambda_{{4}}\lambda_{{12}}-2\,\lambda_{{15}},\lambda_{{5}}-\lambda_{{15}},-\lambda_{{16}}-\lambda_{{17}},-3\,\lambda_{{12}},-2\,\lambda_{{13}}-\lambda_{{19}},0,0,-\lambda_{{21}}\lambda_{{3}}+\lambda_{{22}},\lambda_{{23}}],\\
T_4=&[0,0,0,0,\frac{2}{3}\,\lambda_{{4}},0,0,0,-\frac{1}{3}\,\lambda_{{6}},-\frac{1}{3}\,\lambda_{{7}}-\frac{2}{9}\,{\lambda_{{4}}}^{2},-\frac{1}{3}\,\lambda_{{8}}+\frac{4}{3}\,\lambda_{{9}},-2,0,0,\frac{2}{3}\,\lambda_{{4}},0,\\
&-\frac{1}{3}\,\lambda_{{14}},-\frac{2}{9}\,\lambda_{{4}}\lambda_{{12}}-\frac{1}{3}\,\lambda_{{5}}-\frac{1}{3}\,\lambda_{{15}},0,\frac{5}{3}\,\lambda_{{12}},0,0,0,-\frac{1}{3}\,\lambda_{{21}}\lambda_{{3}}+\frac{1}{3}\,\lambda_{{22}}],\\
T_5=&[0,0,0,0,0,0,0,0,0,0,0,0,-1,0,0,0,0,0,-3,0,0,0,0,0],\\
T_6=&[0,0,0,0,0,0,0,0,0,0,0,0,0,0,0,0,0,0,0,-1,0,0,0,0].\\
\end{array}
\end{equation}
\end{small}
 The solution of system (\ref{c9}) is found in two steps using Maple through the chain (\ref{cad}). First we
consider the following subsystem of equations (\ref{c9})
\begin{equation}\label{c10}
T_i J=0, i=3\dots6.
\end{equation}
In 24-dimensional space of variables $\lambda_i, i=1\dots24$, the
rank of the system (\ref{c10}) is 4, so it has 20 functionally
independent solutions which are given as:
\begin{footnotesize}
\begin{equation}\label{c11}
\begin{array}{ll}
\gamma_{{1}}=&\lambda_{{1}},\\
\gamma_{{2}}=&\lambda_{{2}},\\
\gamma_{{3}}=&\lambda_{{3}} ,\\
\gamma_{{4}}=&\lambda_{{4}},\\
\gamma_{{5}}=&\lambda_{{6}},\\
\gamma_{{6}}=&-{\frac{2\,{\lambda_{{4}}}^{2}\lambda_{{7}}+3\,{\lambda_{{7}}}^{2}-3\,\lambda_{{6}}\lambda_{{8}}}{3\lambda_{{6}}}},\\
\gamma_{{7}}=&-{\frac{2\,{\lambda_{{4}}}^{3}\lambda_{{7}}+3\,{\lambda_{{7}}}^{2}\lambda_{{4}}-3\,{\lambda_{{6}}}^{2}\lambda_{{5}}-6\,\lambda_{{6}}\lambda_{{4}}\lambda_{{9}}}{6\lambda_{{6}}\lambda_{{4}}}} ,\\
\gamma_{{8}}=&{\frac{2\,\lambda_{{6}}\lambda_{{4}}\lambda_{{5}}+6\,\lambda_{{6}}\lambda_{{10}}-6\,\lambda_{{7}}\lambda_{{9}}+3\,\lambda_{{7}}\lambda_{{8}}}{6\lambda_{{6}}}} ,\\
\gamma_{{9}}=&{\frac{4\,\lambda_{{6}}{\lambda_{{4}}}^{3}\lambda_{{7}}\lambda_{{5}}+6\,{\lambda_{{4}}}^{2}{\lambda_{{6}}}^{2}\lambda_{{11}}+12\,{\lambda_{{4}}}^{2}\lambda_{{6}}\lambda_{{7}}\lambda_{{10}}-6\,{\lambda_{{4}}}^{2}{\lambda_{{7}}}^{2}\lambda_{{9}}+3\,{\lambda_{{4}}}^{2}{\lambda_{{7}}}^{2}\lambda_{{8}}+3\,{\lambda_{{6}}}^{2}\lambda_{{4}}\lambda_{{5}}\lambda_{{8}}-12\,{\lambda_{{6}}}^{2}\lambda_{{4}}\lambda_{{5}}\lambda_{{9}}+3\,\lambda_{{6}}\lambda_{{4}}\lambda_{{5}}{\lambda_{{7}}}^{2}-3\,{\lambda_{{6}}}^{3}{\lambda_{{5}}}^{2}}{6{\lambda_{{6}}}^{2}{\lambda_{{4}}}^{2}}} ,\\
\gamma_{{10}}=&{\frac{3\,\lambda_{{5}}+\lambda_{{4}}\lambda_{{12}}}{ \lambda_{{4}}}},\\
\gamma_{{11}}=&\lambda_{{14}},\\
\gamma_{{12}}=&-{\frac {\lambda_{{14}}\lambda_{{7}}+\lambda_{{6}}\lambda_{{5}}-\lambda_{{6}}\lambda_{{15}}}{\lambda_{{6}}}} ,\\
\gamma_{{13}}=&-{\frac{-3\,\lambda_{{16}}{\lambda_{{6}}}^{2}+6\,\lambda_{{6}}\lambda_{{7}}\lambda_{{15}}+2\,\lambda_{{6}}\lambda_{{7}}\lambda_{{4}}\lambda_{{12}}-3\,\lambda_{{14}}{\lambda_{{7}}}^{2}}{3{\lambda_{{6}}}^{2}}} ,\\
\gamma_{{14}}=&{\frac{2\,\lambda_{{4}}{\lambda_{{6}}}^{2}\lambda_{{17}}+2\,\lambda_{{4}}\lambda_{{6}}\lambda_{{7}}\lambda_{{5}}-2\,\lambda_{{4}}\lambda_{{6}}\lambda_{{7}}\lambda_{{15}}+\lambda_{{4}}\lambda_{{14}}{\lambda_{{7}}}^{2}+{\lambda_{{6}}}^{2}\lambda_{{5}}\lambda_{{14}}}{2{\lambda_{{6}}}^{2}\lambda_{{4}}}}\\
\gamma_{{15}}=&{\frac{2\,{\lambda_{{6}}}^{3}\lambda_{{4}}\lambda_{{5}}\lambda_{{12}}+6\,{\lambda_{{6}}}^{3}\lambda_{{4}}\lambda_{{18}}+3\,{\lambda_{{6}}}^{3}{\lambda_{{5}}}^{2}+3\,{\lambda_{{6}}}^{3}\lambda_{{5}}\lambda_{{15}}}{6\lambda_{{4}}{\lambda_{{6}}}^{3}}} ,\\
&+{\frac{-6\,\lambda_{{7}}{\lambda_{{6}}}^{2}\lambda_{{4}}\lambda_{{17}}-6\,\lambda_{{7}}{\lambda_{{6}}}^{2}\lambda_{{4}}\lambda_{{16}}-3\,\lambda_{{7}}{\lambda_{{6}}}^{2}\lambda_{{5}}\lambda_{{14}}-3\,\lambda_{{6}}\lambda_{{4}}\lambda_{{5}}{\lambda_{{7}}}^{2}+9\,{\lambda_{{7}}}^{2}\lambda_{{4}}\lambda_{{6}}\lambda_{{15}}+2\,{\lambda_{{7}}}^{2}{\lambda_{{4}}}^{2}\lambda_{{6}}\lambda_{{12}}-3\,\lambda_{{14}}{\lambda_{{7}}}^{3}\lambda_{{4}}}{6\lambda_{{4}}{\lambda_{{6}}}^{3}}} ,\\
\gamma_{{16}}=&-3\,\lambda_{{13}}+\lambda_{{19}},\\
\gamma_{{17}}=&\lambda_{{21}},\\
\gamma_{{18}}=&\lambda_{{22}},\\
\gamma_{{19}}=&-{\frac {\lambda_{{7}}\lambda_{{21}}\lambda_{{3}}-\lambda_{{7}}\lambda_{{22}}-\lambda_{{23}}\lambda_{{6}}}{\lambda_{{6}}}} ,\\
\gamma_{{20}}=&-{\frac {-2\,\lambda_{{4}}{\lambda_{{6}}}^{2}\lambda_{{24}}-2\,\lambda_{{4}}\lambda_{{7}}\lambda_{{23}}\lambda_{{6}}+\lambda_{{4}}{\lambda_{{7}}}^{2}\lambda_{{21}}\lambda_{{3}}-\lambda_{{4}}{\lambda_{{7}}}^{2}\lambda_{{22}}-{\lambda_{{6}}}^{2}\lambda_{{5}}\lambda_{{21}}\lambda_{{3}}+{\lambda_{{6}}}^{2}\lambda_{{5}}\lambda_{{22}}}{2{\lambda_{{6}}}^{2}\lambda_{{4}}}} ,\\
\end{array}
\end{equation}
\end{footnotesize}
In variables $\gamma_i, i=1\dots20$, the remaining non-zero
operators in system (\ref{c9}) take the form
\begin{small}
\begin{equation}\label{c12}
\begin{array}{ll}
X_1=&[0,1,0,0,0,0,0,0,0,0,0,0,0,0,0,0,0,0,0,0],\\
X_2=&[0,0,1,0,0,0,0,0,0,0,0,0,0,0,0,0,0,\gamma_{{17}},0,0],\\
X_3=&[0,0,\gamma_{{3}},-\gamma_{{4}},-2\,\gamma_{{5}},-2\,\gamma_{{6}},-2\,\gamma_{{7}},-2\,\gamma_{{8}},-2\,\gamma_{{9}},0,-\gamma_{{11}},-\gamma_{{12}},-\gamma_{{13}},-\gamma_{{14}},-\gamma_{{15}},0,0,\gamma_{{18}},\gamma_{{19}},\gamma_{{20}}],\\
T_1=&[1,0,0,0,0,0,0,0,0,0,0,0,0,0,0,0,0,0,0,0],\\
T_2=&,[0,0,-\gamma_{{3}},\gamma_{{4}},3\,\gamma_{{5}},\gamma_{{6}},\gamma_{{7}},0,-\gamma_{{9}},-2\,\gamma_{{10}},0,-\gamma_{{12}},-2\,\gamma_{{13}},-2\,\gamma_{{14}},-3\,\gamma_{{15}},-3\,\gamma_{{16}},\gamma_{{17}},0,-\gamma_{{19}},-2\,\gamma_{{20}}]\\
\end{array}
\end{equation}
\end{small}
Finally, we consider the following subsystem of equations
(\ref{c3})
\begin{equation}\label{c13}
X_i J=0, i=1\dots3, T_k J=0, k=1\dots2.
\end{equation}
In 20-dimensional space of variables $\gamma_i, i=1\dots20$, the
rank of the system (\ref{c13}) is 5, so it has 15 functionally
independent solutions which are given as:
\begin{equation}\label{c14}
\begin{array}{llllll}
\beta_{{1}}={\frac {\gamma_{{6}}\gamma_{{5}}}{{\gamma_{{4}}}^{4}}},&\beta_{{2}}={\frac{\gamma_{{7}}\gamma_{{5}}}{{\gamma_{{4}}}^{4}}},&\beta_{{3}}={\frac{\gamma_{{8}}{\gamma_{{5}}}^{2}}{{\gamma_{{4}}}^{6} }} ,&\beta_{{4}}={\frac{\gamma_{{9}}{\gamma_{{5}}}^{3}}{{\gamma_{{4}}}^{8} }} ,&\beta_{{5}}={\frac{{\gamma_{{5}}}^{2}\gamma_{{10}}}{{\gamma_{{4}}}^{4 }}} ,&\beta_{{6}}={\frac{\gamma_{{5}}\gamma_{{11}}}{{\gamma_{{4}}}^{3}}},\\
\beta_{{7}}={\frac{\gamma_{{12}}{\gamma_{{5}}}^{2}}{{\gamma_{{4}}}^{5}}},&\beta_{{8}}={\frac{\gamma_{{13}}{\gamma_{{5}}}^{3}}{{\gamma_{{4}}}^{7}}},&\beta_{{9}}={\frac{\gamma_{{14}}{\gamma_{{5}}}^{3}}{{\gamma_{{4}}}^{7}}},&\beta_{{10}}={\frac{\gamma_{{15}}{\gamma_{{5}}}^{4}}{{\gamma_{{4}}}^{9}}},&\beta_{{11}}={\frac{\gamma_{{16}}{\gamma_{{5}}}^{3}}{{\gamma_{{4}}}^{6}}}\\
\end{array}
\end{equation}
and
\begin{equation}\label{c15}
\begin{array}{llll}
\beta_{{12}}={\frac{\gamma_{{17}}{\gamma_{{4}}}^{2}}{\gamma_{{5}}}},&\beta_{{13}}=-{\frac { \left(\gamma_{{17}}\gamma_{{3}}-\gamma_{{18}} \right) {\gamma_{{4}}}^{3}}{\gamma_{{5}}}},&\beta_{{14}}=\gamma_{{19}}\gamma_{{4}},&\beta_{{15}}={\frac {\gamma_{{5}}\gamma_{{20}}}{\gamma_{{4}}}}.\\
\end{array}
\end{equation}
Here $\beta_{{12}}, \beta_{{13}}, \beta_{{14}}$ and $\beta_{{15}}$
are the only invariants depending on the variables $K, L, M$ and
$N$. Then the general solution of (\ref{ccc}), for
$\xi=\xi(x),\eta=\eta(x,y)$, can be given implicitly by back
substitution as
\begin{equation}\label{c16}
\begin{array}{llll}
{\frac {{f^2_{{4,4}}}}{f_{{4,4,4}}}}K={\it F_1},\\
{\frac {{f^3_{{4,4}}}}{f_{{4,4 ,4}}}} \left( p K_{{}}-L_{{}} \right)={\it F_2} ,\\
{\frac {{f_{{4,4}}}}{f_{{4,4,4}}}}\left((f_4f_{4,4,4}+3f_{3,4,4})(pK-L)+3f_{4,4,4}(qK-M)\right)={\it F_3} ,\\
\frac {1}{f^2_{4,4}f_{{4,4,4}}} ( (f_{{4,4}}f_{{3,4,4}}(2\,f_4 f_{{4,4,4}}+3{f_{{3,4,4}}})-{f_{{4,4,4}}}^{2}( 2f_{{3,4}}f_{{4}}-6\,f_{{2,4}}+3\,f_{{3,3}} ) ) (p K_{{}}-L_{{}} ) \\
+6\,f_{{4,4}}f_{{3,4,4}}f_{{4,4,4}}(qK-M)+6f_{{4,4}}{f_{{4,4,4}}}^{2}(\,f_{{}}K_{{}}-\,N_{{}} ))=F_4.\\
\end{array}
\end{equation}
where $F_1, F_2, F_3$ and $F_4$ are the arbitrary functions of
$\beta_{{i}}, i=1...11.$

Finally, solving system (\ref{c16}) gives the variables $K, L, M$
and $N$ in terms of four arbitrary functions $F_1, F_2, F_3$ and
$F_4$ which provide four independent invariant differentiation
operators $\mathcal{D}_1,\mathcal{D}_2,\mathcal{D}_3$ and
$\mathcal{D}_4$ via (\ref{cc}).
\endproof
\section{Illustrative Examples of Equivalent Equations}
We present illustrative examples of third-order ODEs, not
quadratic in the second-order derivative, taken from the works
\cite{Mahomed1996,Zaitsev2002,Dunajski2013}. In these studies, the
symmetry algebra properties were investigated in
\cite{Mahomed1996}, exact solutions for certain classes of
third-order ODEs in \cite{Zaitsev2002} and a characterization of
Lorentzian three-dimensional hyper-CR Einstein-Weyl structures in
terms of invariants of the associated third-order ODEs in
\cite{Dunajski2013}.
\begin{example}\rm \label{ex1}
The invariants obtained here may be used when we need to prove the
nonequivalence of two given equations under fiber preserving point transformations. Consider the equation
\begin{equation}\label{eee1}
\bar{y}'''=\bar{y}''^4,
\end{equation}
with the invariants, given by Theorem 3.2, as
\begin{equation}\label{eee2}
\begin{array}{l}
\bar{\beta}_1=0,\bar{\beta}_2=0,\bar{\beta}_3=0,\bar{\beta}_4=0,\bar{\beta}_5=\frac{5}{27},\bar{\beta}_6=\frac{5}{9},\bar{\beta}_7=0,\bar{\beta}_8=0,\bar{\beta}_9=0,\bar{\beta}_{10}=0,\bar{\beta}_{11}=\frac{5}{243},\\
\end{array}
\end{equation}
and the equation
\begin{equation}\label{eee3}
y'''={y''}^3,
\end{equation}
with the invariants, via Theorem 3.2,
\begin{equation}\label{eee4}
\begin{array}{l}
\beta_1=0,\beta_2=0,\beta_3=0,\beta_4=0,\beta_5=\frac{1}{12},\beta_6=\frac{1}{3},\beta_7=0,\beta_8=0,\beta_9=0,\beta_{10}=0,\beta_{11}=0.\\
\end{array}
\end{equation}

Since their respective invariants are not equal, these two equations
(\ref{eee1}) and (\ref{eee3}) are obviously nonequivalent with
respect to the fiber preserving transformations $\bar{x}=\phi(x),
\bar{y}=\psi(x,y)$.
\end{example}
\begin{example}\rm  \label{ex2}
Consider the equation \cite[Section 3.2]{Zaitsev2002}, also listed
in \cite[Section 8.3.3]{Mahomed1996},
\begin{equation}\label{eee5}
\bar{y}'''=A~\bar{y}''^\delta,~\delta\ne0,1,2
\end{equation}
and by setting $\delta=3$ this becomes
\begin{equation}
\bar y'''=A~\bar y''^3.\label{2}
\end{equation}
 We know the invariants of (\ref{2}) by Theorem 3.2 which are given in (\ref{eee4}). By means
of the fiber preserving transformation
\begin{equation}
\bar x=\ln x,\; \bar y=x+y\label{3}
\end{equation}
 this equation (\ref{2}) transforms to
\begin{equation}
 x^3y'''+3x^2y''+xy'+x=A(x+xy'+x^2y'')^3.\label{4}
\end{equation}
Equation (\ref{4}) also has the same values of the invariants as
given in (\ref{eee4}).

Now we focus our attention on the third-order ODE
\begin{equation} {3\over y'^5}y''^2-{1\over
y'^4} y'''=-A\left({y''\over y'^3} \right)^3 \label{5}
\end{equation}
 This equation (\ref{5}) maps to (\ref{2}) via the interchange transformation
\begin{equation}
\bar x=y,\; \bar y =x\label{6}
\end{equation}
 but this transformation is not fiber preserving. The invariants for
this ODE (\ref{5}) are not identically constant and so there is no
fiber preserving transformation between ODEs (\ref{5}) and
(\ref{2}).
\end{example}
\begin{remark}\rm
The special case, of the third order ODE (\ref{eee5}),
\begin{equation}
\bar y'''=\bar y''^{\frac{3}{2}}\label{7}
\end{equation}
defines Einstein-Weyl geometry of hyper-CR type and is of recent
interest in physics \cite{Dunajski2013}. It can be characterized
by the invariants, by Theorem 3.2, as
\begin{equation}\label{8}
\begin{array}{l}
\bar{\beta}_i=0,~i=1..11.
\end{array}
\end{equation}
\end{remark}
\begin{example}\rm  \label{ex3}
Consider the equation \cite[Section 3.2]{Zaitsev2002}
\begin{equation}\label{eee7}
\bar{y}'''=A~\bar{x}^\alpha~\bar{y}''^\delta,~\delta\ne0,1,2
\end{equation}
with the only nonzero third order invariants, by Theorem 3.2, given by
\begin{equation}\label{eee8}
\begin{array}{l}
\bar{\beta}_5=\,{\frac{\left(2\,\delta-3 \right)  \left( \delta-2 \right)^{2}}{3\delta\,\left(\delta-1 \right)^{2}}}+{\frac {\alpha\, \left( \delta-2 \right) ^{2}}{A\delta\left( \delta-1\right) ^{2}}}~s,\bar{\beta}_6={\frac {2 \left( \delta-2 \right)\left(2 \delta-3\right) }{3 \delta\left( \delta-1 \right)}}+{\frac {\alpha \left( \delta-2 \right) }{A \delta\left(\delta-1 \right)}}~s,\\
\bar{\beta}_{11}={\frac{2\left( \delta-2 \right) ^{3} \left(\delta-3 \right)\left(2 \delta-3\right)}{9{\delta}^{2} \left(\delta-1 \right) ^{3}}}+{\frac {\alpha\, \left( \delta-2 \right)^{3}\left( \delta-3 \right)}{A\delta^2\left( \delta-1\right)^{3}}}~s+{\frac {\alpha\, \left(\alpha-1 \right)  \left( \delta-2 \right) ^{3}}{A^2{\delta}^{2} \left( \delta-1 \right)^{3}}}~s^2,\\
\end{array}
\end{equation}
where $s=\bar{x}^{(-1-\alpha)}\bar{y}''^{(1-\delta)}$.

It is clear here that the invariants (\ref{eee8}) are not
identically constant. However, their image using a fiber
preserving transformation should match the corresponding
differential invariants of the transformed equation using the same
fiber preserving transformation such that
\begin{equation}\label{eee9}
\begin{array}{l}
\bar{\beta}_5(\bar x,\bar y,\bar y',\bar y'')=\beta_5(x,y,y',y''),~\bar{\beta}_6(\bar x,\bar y,\bar y',\bar y'')=\beta_6(x,y,y',y''),~\bar{\beta}_{11}(\bar x,\bar y,\bar y',\bar y'')=\beta_{11}(x,y,y',y'').\\
\end{array}
\end{equation}

We take $\alpha=3,~\delta=4$ for illustration. Then write this as
the transformed third-order ODE
\begin{equation}
\bar y'''=A~\bar x ^3 \bar y''^4.\label{eee10}
\end{equation}
 The only nonzero third order invariants of (\ref{eee10}), by Theorem 3.2,
 are the following polynomials
\begin{equation}\label{eee11}
\begin{array}{l}
\bar{\beta}_5=\frac{5}{27}+\frac{1}{3A}~s,\bar{\beta}_6=\frac{5}{9}+\frac{1}{2A}~s,\bar{\beta}_{11}=\frac{5}{243}+\frac{1}{18A}~s+\frac{1}{9A^2}~s^2,\\
\end{array}
\end{equation}
where $s=\bar{x}^{-4}\bar{y}''^{-3}$.

By means of the fiber preserving transformation
\begin{equation}
\bar x= \frac{1}{x},\; \bar y=\frac{y}{x}\label{eee12}
\end{equation}
 this equation (\ref{eee10}) transforms to
\begin{equation}
y'''=-A~{x}^{4}{y''}^{4}-3\,\frac{y''}{x}.\label{eee13}
\end{equation}
The only nonzero third order invariants of (\ref{eee13}), by
Theorem 3.2, are
\begin{equation}\label{eee14}
\begin{array}{l}
\beta_5=\frac{5}{27}+\frac{1}{3A}~\bar s,\beta_6=\frac{5}{9}+\frac{1}{2A}~\bar s,\beta_{11}=\frac{5}{243}+\frac{1}{18A}~\bar s+\frac{1}{9A^2}~\bar s^2,\\
\end{array}
\end{equation}
where $\bar s=x^{-5} y''^{-3}$.

Since the transformation (\ref{eee12}) transforms the variable $s$
to the variable $\bar s$, then the system (\ref{eee9}) is
satisfied.
\end{example}
\begin{remark}\rm
Using (\ref{eee11}), it is clear that the variable $s$ is
invariant and can be given as $s=A(3\bar \beta_5-\frac{5}{9})$.
Therefore, one can use any symbolic package such as Maple or
Mathematica to study the equivalence of any ODE
$y'''=f(x,y,y',y'')$ to the ODE (\ref{eee10}) via the fiber
preserving transformations $\bar{x}=\phi(x),~\bar{y}=\psi(x,y)$.
This can be done by comparing the differential invariants
(\ref{d4}) calculated directly from the ODE $y'''=f(x,y,y',y'')$
and the invariants (\ref{eee11}) after replacing the variables $s$
by $\bar s=A(3 \beta_5-\frac{5}{9})$. Moreover, the relation
$s=\bar s$ may help in constructing the fiber preserving
transformation.
\end{remark}
\begin{remark}\rm
The invariants (\ref{eee11}) are not identically constants.
However, one can use them to construct identically constant
invariants by eliminating the variable $s$. For examples, the
following invariants $\bar I_1$ and $\bar I_2$ are identically
zero
\begin{equation}
\begin{array}{l}
\bar I_1=3\bar \beta_5-2\bar \beta_6+\frac{5}{9}=0,\bar I_2=\frac{5}{243}+\frac{1}{18}~(3\bar \beta_5-\frac{5}{9})+\frac{1}{9}~(3\bar \beta_5-\frac{5}{9})^2-\bar{\beta}_{11}=0.\\
\end{array}
\end{equation}
\end{remark}
\begin{example}\rm  \label{ex4}
Consider the equation \cite[Section 3.2]{Zaitsev2002}
\begin{equation}\label{eee15}
\bar{y}'''=A~\bar{y}^{(-7)}~\bar{y}'^7~\bar{y}''^3.
\end{equation}
The third order invariants, by Theorem 3.2, are the following
polynomials
\begin{equation}\label{eee16}
\begin{array}{l}
\bar{\beta}_1=-{\frac {7}{9}}\,{s}^{2}-{\frac {7}{18}}\,s,\\
\bar{\beta}_2=-{\frac {7}{72}}\,{s}^{3}t-{\frac{7}{12}}\,{s}^{2}-{\frac {7}{24}}\,s,\\
\bar{\beta}_3={\frac {49}{432}} \,{s}^{4}t+{\frac{7}{27}}\,{s}^{3}+{\frac {35}{432}}\,{s}^{2},\\
\bar{\beta}_4={\frac {35}{2592}}\,{s}^{6}{t}^{2}+{\frac{245}{1296}}\,t{s}^{5}+\left( {\frac {245}{648}}+{\frac {49}{1944}}\,t \right) {s}^{4}+{\frac {161}{1296}}\,{s}^{3}+{\frac {7}{864}}\,{s}^{2},\\
\bar{\beta}_5={\frac {7}{12}}\,{s}^{3}t+ \left( {\frac{7}{12}}-{\frac {7}{12}} \,t \right) {s}^{2}+{\frac {7}{12}}\,s+\frac{1}{12},\\
\bar{\beta}_6=-\frac{7}{6}\,t{s}^{2}+\frac{7}{6}\,s+\frac{1}{3},\\
\bar{\beta}_7={\frac{7}{12}}\,{s}^{3} t-{\frac {35}{36}}\,{s}^{2}-{\frac {7}{36}}\,s,\\
\bar{\beta}_8=-{\frac {245}{648}}\,{s}^{4}t+ \left( {\frac{7}{18}}+{\frac {49}{216} }\,t \right) {s}^{3}-{\frac{35}{216}}\,{s}^{2}-{\frac {7}{216}}\,s,\\
\bar{\beta}_9={\frac {7}{432}}\,{t}^{2}{s}^{5}+{\frac{49}{1296}}\,{s}^{4}t+ \left( {\frac {49}{108}}+{\frac{7}{432}}\,t \right) {s}^{3}+{\frac {35}{432} }\,{s}^{2},\\
\bar{\beta}_{10}=-{\frac {217}{7776}}\,{s}^{6}{t}^{2}-{\frac{343}{2592}}\,t{s}^{5}+\left( -{\frac {49}{432}}-{\frac {245}{3888}}\,t \right) {s}^{4}+{\frac {175}{2592}}\,{s}^{3}+{\frac {35}{2592}}\,{s}^{2},\\
\bar{\beta}_{11}={\frac{7}{9}}\,{t}^{2}{s}^{4}-{\frac {35}{36}}\,{s}^{3}t,\\
\end{array}
\end{equation}
where $s=\bar y^{7}\bar y'^{(-8)} \bar y''^{(-1)}$ and $t=\bar
y^{(-8)}\bar y'^{(10)}$.

By means of the fiber preserving transformation
\begin{equation}
\bar x= \frac{x}{x-1},\; \bar y=\frac{y}{x-1}\label{eee17}
\end{equation}
 this equation (\ref{eee15}) transforms to
\begin{equation}
y'''={\frac { y''\left(  \left( y'  x-y' -y \right) ^{7} \left( x-1 \right)^{12}{y''}^{2}-3\,{y}^{7}\right)}{{y}^{7} \left( x-1 \right)}}.\label{eee18}\\
\end{equation}
The third order invariants calculated directly from the ODE
(\ref{eee18}), by Theorem 3.2, match the invariants
(\ref{eee16}) after replacing the variables $s$ and $t$ by
their image $\bar s=(x-1)^{(-10)} y^{7}(y'-x
y'+y)^{(-8)}y''^{(-1)}$ and $\bar t=(x-1)^{8} y^{(-8)}(y'-x
y'+y)^{10}$ under the transformation (\ref{eee17}).
\end{example}
This completes the examples.

\section{Conclusion}
We have provided an extension of the work of
Bagderina \cite{Bagderina2008} who solved the equivalence problem for
scalar third-order ordinary differential equations (ODEs), quadratic
in the second-order derivative, via point transformations. Here
we considered the equivalence problem for third-order ODEs of the general form
$y'''=f(x,y,y',y'')$, which are not quadratic in the second-order
derivative, under the pseudo-group of fiber preserving equivalence transformations
$\bar{x}=\phi(x), \bar{y}=\psi(x,y)$. We utilized Lie's infinitesimal method to
obtain the differential invariants of this general class of ODEs with $f$ not quadratic in $y''$. All third
order differential invariants of this pseudo-group and the invariant
differentiation operators were determined. These are stated as two Theorems 3.1 and 3.2
in Section 3. These yield simple
necessary explicit conditions for a third-order ODE to be equivalent to
the respective canonical form under pseudo-group of point transformations.
As illustrative examples, we gave
a number of equations from the existing literature in Section 4.

It would be of further interest to look at the equivalence problem
under a more general equivalence group of point transformations.
This necessitate further coding in Maple or Mathematica in order
to facilitate the large calculations which also have prevailed for
the pseudo-group of fiber preserving point transformations
considered in this work.
\subsection*{Acknowledgments}
 The authors would like to thank the King Fahd University of
Petroleum and Minerals for its support and excellent research
facilities. FM thanks the NRF of South Africa for  research support through a grant.
\section*{Appendix A: The differential operators of the homogeneous linear system of
PDEs (\ref{c3}) }
\begin{footnotesize}
\begin{eqnarray*}
X_1 &=& [0,1,0,0,0,0,0,0,0,0,0,0,0,0,0,0,0,0,0,0,0,0,0,0,0,0,0,0,0,0,0,0,0,0,0,0,0, \\
  && 0,0,0,0,0,0]  \\
X_2  &=& [0,0,1,0,0,0,0,0,0,0,0,0,0,0,0,0,0,0,0,0,0,0,0,0,0,0,0,0,0,0,0,0,0,0,0,0,0, \\
  && 0,0,0,z_{{40}},0,0] \\
X_3 &=& [0,0,z_{{3}},z_{{4}},z_{{5}},0,0,0,-z_{{9}},-z_{{10}},-z_{{11}},-z_{{12}},-z_{{13}},-z_{{14}},-2\,z_{{15}},-2\,z_{{16}},-2\,z_{{17}},-2\,z_{{18}}, \\
 && -2\,z_{{19}},-2\,z_{{20}},-2\,z_{{21}},-2\,z_{{22}},-2\,z_{{23}},-2\,z_{{24}},z_{{25}},0,0,0,-z_{{29}},-z_{{30}},-z_{{31}},-z_{{32}},-z_{{33}},\\
 && -z_{{34}},z_{{35}},0,0,0,z_{{39}},0,z_{{41}},z_{{42}},z_{{43}}] \\
X_4  &=& [0,0,0,1,0,0,0,0,0,0,0,0,0,0,0,0,0,0,0,0,0,0,0,0,0,0,0,0,0,0,0,0,0,0,0,0,0,\\
&& 0,0,0,0,z_{{40}},0]\\
X_5 &=& [0,0,0,2\,z_{{3}},3\,z_{{4}},-z_{{7}},-2\,z_{{8}},3,-2\,z_{{10}},-2\,z_{{11}}-z_{{12}},-z_{{13}},-4\,z_{{13}},-2\,z_{{14}},0,-3\,z_{{16}}, -2\,z_{{18}}\\
&&-2\,z_{{17}},-2\,z_{{19}},-z_{{21}}-4\,z_{{19}},-2\,z_{{20}}-z_{{22}},-z_{{23}},-6\,z_{{22}},-4\,z_{{23}},-2\,z_{{24}},0,4\,z_{{5}},-z_{{27}},-2\,z_{{28}},\\
&&0,-z_{{9}}-2\,z_{{30}},-z_{{32}}-z_{{10}}-2\,z_{{31}},-z_{{11}}-z_{{33}},-z_{{12}}-4\,z_{{33}},-z_{{13}}-2\,z_{{34}},-z_{{14}},5\,z_{{25}},-z_{{37}},\\
&&-2\,z_{{38}},0,6\,z_{{35}},0,0,z_{{41}}+z_{{40}}z_{{3}},2\,z_{{42}}+z_{{40}}z_{{4}}]  \\
X_6  &=& [0,0,0,{z_{{3}}}^{2},3\,z_{{3}}z_{{4}},z_{{5}}-z_{{8}}z_{{4}}-z_{{7}}z_{{3}},3\,z_{{4}}-2\,z_{{8}}z_{{3}},3\,z_{{3}},z_{{6}}-2\,z_{{3}}z_{{10}}-2\,z_{{4}}z_{{11}}, -z_{{4}}z_{{13}}-2\,z_{{3}}z_{{11}}\\
&&-z_{{3}}z_{{12}},-z_{{3}}z_{{13}}-z_{{4}}z_{{14}},-2\,z_{{8}}-4\,z_{{3}}z_{{13}},3-2\,z_{{3}}z_{{14}},0,-3\,z_{{3}}z_{{16}}-3\,z_{{4}}z_{{17}},-2\,z_{{4}}z_{{19}}-2\,z_{{3}}z_{{17}}\\
&&-2\,z_{{3}}z_{{18}}-z_{{10}},-z_{{11}}-2\,z_{{3}}z_{{19}}-2\,z_{{4}}z_{{20}}, -z_{{12}}-z_{{3}}z_{{21}}-2\,z_{{11}}-4\,z_{{3}}z_{{19}}-z_{{4}}z_{{22}},-z_{{13}}-z_{{3}}z_{{22}}\\
&&-2\,z_{{3}}z_{{20}}-z_{{4}}z_{{23}}, -z_{{4}}z_{{24}}-z_{{14}}-z_{{3}}z_{{23}},-6\,z_{{13}}-6\,z_{{3}}z_{{22}},-2\,z_{{14}}-4\,z_{{3}}z_{{23}},-2\,z_{{3}}z_{{24}},0,3\,{z_{{4}}}^{2}\\
&&+4\,z_{{5}}z_{{3}}, z_{{25}}-z_{{5}}z_{{8}}-z_{{4}}z_{{28}}-z_{{7}}z_{{4}}-z_{{3}}z_{{27}},3\,z_{{5}}-2\,z_{{3}}z_{{28}}-2\,z_{{8}}z_{{4}},3\,z_{{4}},-z_{{3}}z_{{9}}+z_{{26}}-2\,z_{{4}}z_{{31}}\\
&&-2\,z_{{3}}z_{{30}}-2\,z_{{4}}z_{{10}}-2\,z_{{5}}z_{{11}}, -z_{{5}}z_{{13}}-z_{{4}}z_{{12}}-z_{{3}}z_{{10}}-2\,z_{{3}}z_{{31}}-2\,z_{{4}}z_{{11}}-z_{{3}}z_{{32}}-z_{{4}}z_{{33}},-z_{{4}}z_{{13}}\\
&&-z_{{3}}z_{{11}}-z_{{3}}z_{{33}}-z_{{5}}z_{{14}}-z_{{4}}z_{{34}},-2\,z_{{28}}-4\,z_{{4}}z_{{13}}-z_{{3}}z_{{12}}-4\,z_{{3}}z_{{33}},-2\,z_{{3}}z_{{34}}-z_{{3}}z_{{13}}-2\,z_{{4}}z_{{14}},\\
&&-z_{{3}}z_{{14}},10\,z_{{5}}z_{{4}}+5\,z_{{3}}z_{{25}}, z_{{35}}-2\,z_{{5}}z_{{28}}-2\,z_{{4}}z_{{27}}-z_{{8}}z_{{25}}-z_{{3}}z_{{37}}-z_{{4}}z_{{38}}-z_{{5}}z_{{7}},3\,z_{{25}}-2\,z_{{5}}z_{{8}}\\
&&-4\,z_{{4}}z_{{28}}-2\,z_{{3}}z_{{38}},3\,z_{{5}},10\,{z_{{5}}}^{2}+6\,z_{{3}}z_{{35}}+15\,z_{{4}}z_{{25}},0,0,z_{{41}}z_{{3}},z_{{41}}z_{{4}}+2\,z_{{42}}z_{{3}}]
\end{eqnarray*}
\end{footnotesize}
\begin{footnotesize}
\begin{eqnarray*}
X_7 &=&  [0,0,0,0,1,0,0,0,0,0,0,0,0,0,0,0,0,0,0,0,0,0,0,0,0,0,0,0,0,0,0,0,0,0,0,0,0,\\
&& 0,0,0,0,0,z_{{40}}] \\
X_8  &=& [0,0,0,0,3\,z_{{3}},-z_{{8}},3,0,-2\,z_{{11}},-z_{{13}},-z_{{14}},0,0,0,-3\,z_{{17}},-2\,z_{{19}},-2\,z_{{20}},-z_{{22}},-z_{{23}},\\
&& -z_{{24}},0,0,0,0,6\,z_{{4}},-z_{{28}}-z_{{7}},-2\,z_{{8}},3,-2\,z_{{10}}-2\,z_{{31}},-2\,z_{{11}}-z_{{12}}-z_{{33}},-z_{{13}}-z_{{34}},\\
&& -4\,z_{{13}},-2\,z_{{14}},0,10\,z_{{5}},-2\,z_{{27}}-z_{{38}},-4\,z_{{28}},0,15\,z_{{25}},0,0,0,z_{{41}}+2\,z_{{40}}z_{{3}}] \\
X_9&=&[0,0,0,0,3\,{z_{{3}}}^{2},3\,z_{{4}}-2\,z_{{8}}z_{{3}},6\,z_{{3}},0,-4\,z_{{3}}z_{{11}}-z_{{7}},-2\,z_{{3}}z_{{13}}-2\,z_{{8}},3-2\,z_{{3}}z_{{14}},6,0,\\
&&0,-6\,z_{{3}}z_{{17}}-3\,z_{{10}},-4\,z_{{3}}z_{{19}}-z_{{12}}-4\,z_{{11}},-4\,z_{{3}}z_{{20}}-z_{{13}},-2\,z_{{3}}z_{{22}}-4\,z_{{13}},-2\,z_{{3}}z_{{23}}\\
&&-2\,z_{{14}},-2\,z_{{3}}z_{{24}},0,0,0,0,12\,z_{{3}}z_{{4}},-3\,z_{{8}}z_{{4}}-2\,z_{{7}}z_{{3}}-2\,z_{{3}}z_{{28}}+4\,z_{{5}}, -4\,z_{{8}}z_{{3}}+9\,z_{{4}},\\
&&6\,z_{{3}},-6\,z_{{4}}z_{{11}}-4\,z_{{3}}z_{{10}}-4\,z_{{3}}z_{{31}}-z_{{27}}+z_{{6}},-3\,z_{{4}}z_{{13}}-4\,z_{{3}}z_{{11}}-2\,z_{{3}}z_{{12}}-2\,z_{{3}}z_{{33}}\\
&&-2\,z_{{28}},-3\,z_{{4}}z_{{14}}-2\,z_{{3}}z_{{13}}-2\,z_{{3}}z_{{34}},-8\,z_{{3}}z_{{13}}-2\,z_{{8}},-4\,z_{{3}}z_{{14}}+3,0,20\,z_{{5}}z_{{3}}+15\,{z_{{4}}}^{2},\\
&&-2\,z_{{3}}z_{{38}}-3\,z_{{7}}z_{{4}}-4\,z_{{5}}z_{{8}}-4\,z_{{3}}z_{{27}}-6\,z_{{4}}z_{{28}}+5\,z_{{25}},-6\,z_{{8}}z_{{4}}-8\,z_{{3}}z_{{28}}+12\,z_{{5}},9\,z_{{4}},\\
&&30\,z_{{3}}z_{{25}}+60\,z_{{5}}z_{{4}}, 0,0,0,z_{{40}}{z_{{3}}}^{2}+2\,z_{{41}}z_{{3}}] \\
X_{10}  &=& [0,0,0,0,{z_{{3}}}^{3},-z_{{8}}{z_{{3}}}^{2}+3\,z_{{3}}z_{{4}},3\,{z_{{3}}}^{2},0,-z_{{8}}z_{{4}}-2\,{z_{{3}}}^{2}z_{{11}}-z_{{7}}z_{{3}}+z_{{5}},-{z_{{3}}}^{2}z_{{13}}-2\,z_{{8}}z_{{3}}\\
&&+3\,z_{{4}},3\,z_{{3}}-{z_{{3}}}^{2}z_{{14}},6\,z_{{3}},0,0,-3\,z_{{3}}z_{{10}}-3\,z_{{4}}z_{{11}}-3\,{z_{{3}}}^{2}z_{{17}}+2\,z_{{6}},-2\,{z_{{3}}}^{2}z_{{19}}-4\,z_{{3}}z_{{11}}\\
&&-z_{{4}}z_{{13}}-z_{{3}}z_{{12}},-2\,{z_{{3}}}^{2}z_{{20}}-z_{{3}}z_{{13}}-z_{{4}}z_{{14}}, -4\,z_{{3}}z_{{13}}-{z_{{3}}}^{2}z_{{22}}-2\,z_{{8}},3-2\,z_{{3}}z_{{14}}-{z_{{3}}}^{2}z_{{23}},\\
&&-{z_{{3}}}^{2}z_{{24}},6,0,0,0,6\,{z_{{3}}}^{2}z_{{4}},-{z_{{3}}}^{2}z_{{28}}+4\,z_{{5}}z_{{3}}-3\,z_{{8}}z_{{3}}z_{{4}}-z_{{7}}{z_{{3}}}^{2}+3\,{z_{{4}}}^{2},-2\,z_{{8}}{z_{{3}}}^{2}+9\,z_{{3}}z_{{4}},\\
&&3\,{z_{{3}}}^{2}, -2\,{z_{{3}}}^{2}z_{{31}}-2\,{z_{{3}}}^{2}z_{{10}}+z_{{6}}z_{{3}}-6\,z_{{3}}z_{{4}}z_{{11}}-z_{{3}}z_{{27}}-z_{{4}}z_{{28}}-z_{{5}}z_{{8}}-z_{{7}}z_{{4}}+z_{{25}},-2\,z_{{8}}z_{{4}}\\
&&-2\,{z_{{3}}}^{2}z_{{11}}-2\,z_{{3}}z_{{28}}-3\,z_{{3}}z_{{4}}z_{{13}}-{z_{{3}}}^{2}z_{{33}}-{z_{{3}}}^{2}z_{{12}}+3\,z_{{5}}, -{z_{{3}}}^{2}z_{{13}}-3\,z_{{3}}z_{{4}}z_{{14}}-{z_{{3}}}^{2}z_{{34}}+3\,z_{{4}},\\
&&-4\,{z_{{3}}}^{2}z_{{13}}-2\,z_{{8}}z_{{3}}+6\,z_{{4}},3\,z_{{3}}-2\,{z_{{3}}}^{2}z_{{14}}, 0,15\,z_{{3}}{z_{{4}}}^{2}+10\,z_{{5}}{z_{{3}}}^{2},-3\,z_{{8}}{z_{{4}}}^{2}-2\,{z_{{3}}}^{2}z_{{27}}+10\,z_{{5}}z_{{4}}\\
&&+5\,z_{{3}}z_{{25}}-4\,z_{{5}}z_{{8}}z_{{3}}-3\,z_{{7}}z_{{3}}z_{{4}} -6\,z_{{3}}z_{{4}}z_{{28}}-{z_{{3}}}^{2}z_{{38}},-4\,{z_{{3}}}^{2}z_{{28}}+12\,z_{{5}}z_{{3}}-6\,z_{{8}}z_{{3}}z_{{4}}+9\,{z_{{4}}}^{2},\\
&&9\,z_{{3}}z_{{4}},15\,{z_{{3}}}^{2}z_{{25}}+60\,z_{{5}}z_{{3}}z_{{4}}+15\,{z_{{4}}}^{3},0,0,0,z_{{41}}{z_{{3}}}^{2}]
\end{eqnarray*}
\end{footnotesize}
\begin{footnotesize}
\begin{eqnarray*}
X_{11} &=& [0,0,0,0,0,0,0,0,0,0,0,0,0,0,0,0,0,0,0,0,0,0,0,0,1,0,0,0,0,0,0,0,0,0,0,0,0,\\
&& 0,0,0,0,0,0]  \\
X_{12}  &=& [0,0,0,0,0,1,0,0,0,0,0,0,0,0,0,0,0,0,0,0,0,0,0,0,4\,z_{{3}},-z_{{8}},3,0,-2\,z_{{11}},-z_{{13}},\\
&& -z_{{14}},0,0,0,10\,z_{{4}},-z_{{7}}-2\,z_{{28}},-2\,z_{{8}},3,20\,z_{{5}},0,0,0,0] \\
X_{13}  &=& [0,0,0,0,0,3\,z_{{3}},0,0,-z_{{8}},3,0,0,0,0,-3\,z_{{11}},-z_{{13}},-z_{{14}},0,0,0,0,0,0,0,6\,{z_{{3}}}^{2}, 6\,z_{{4}}\\
&&-3\,z_{{8}}z_{{3}},9\,z_{{3}},0,-z_{{7}}-z_{{28}}-6\,z_{{3}}z_{{11}},-2\,z_{{8}}-3\,z_{{3}}z_{{13}},3-3\,z_{{3}}z_{{14}},6,0,0,30\,z_{{3}}z_{{4}},10\,z_{{5}}\\
&&-6\,z_{{3}}z_{{28}}-3\,z_{{7}}z_{{3}}-6\,z_{{8}}z_{{4}},18\,z_{{4}}-6\,z_{{8}}z_{{3}},9\,z_{{3}},45\,{z_{{4}}}^{2}+60\,z_{{5}}z_{{3}},0,0,0,0] \\
X_{14} &=&  [0,0,0,0,0,3\,{z_{{3}}}^{2},0,0,3\,z_{{4}}-2\,z_{{8}}z_{{3}},6\,z_{{3}},0,0,0,0,-z_{{7}}-6\,z_{{3}}z_{{11}},-2\,z_{{3}}z_{{13}}-2\,z_{{8}},\\
&& 3-2\,z_{{3}}z_{{14}},6,0,0,0,0,0,0,4\,{z_{{3}}}^{3},12\,z_{{3}}z_{{4}}-3\,z_{{8}}{z_{{3}}}^{2},9\,{z_{{3}}}^{2},0, 4\,z_{{5}}-2\,z_{{7}}z_{{3}}-3\,z_{{8}}z_{{4}}\\
&&-6\,{z_{{3}}}^{2}z_{{11}}-2\,z_{{3}}z_{{28}},9\,z_{{4}}-4\,z_{{8}}z_{{3}}-3\,{z_{{3}}}^{2}z_{{13}}, 6\,z_{{3}}-3\,{z_{{3}}}^{2}z_{{14}},12\,z_{{3}},0,0,30\,{z_{{3}}}^{2}z_{{4}},\\
&&-12\,z_{{8}}z_{{3}}z_{{4}}+15\,{z_{{4}}}^{2}+20\,z_{{5}}z_{{3}}-6\,{z_{{3}}}^{2}z_{{28}}-3\,z_{{7}}{z_{{3}}}^{2},36\,z_{{3}}z_{{4}}-6\,z_{{8}}{z_{{3}}}^{2},9\,{z_{{3}}}^{2},90\,z_{{3}}{z_{{4}}}^{2}\\
&&+60\,z_{{5}}{z_{{3}}}^{2},0,0,0,0] \\
X_{15}  &=& [0,0,0,0,0,{z_{{3}}}^{3},0,0,-z_{{8}}{z_{{3}}}^{2}+3\,z_{{3}}z_{{4}},3\,{z_{{3}}}^{2},0,0,0,0,-z_{{8}}z_{{4}}-z_{{7}}z_{{3}}-3\,{z_{{3}}}^{2}z_{{11}}+z_{{5}},\\
&& -{z_{{3}}}^{2}z_{{13}}-2\,z_{{8}}z_{{3}}+3\,z_{{4}},3\,z_{{3}}-{z_{{3}}}^{2}z_{{14}},6\,z_{{3}},0,0,0,0,0,0,{z_{{3}}}^{4},6\,{z_{{3}}}^{2}z_{{4}}-z_{{8}}{z_{{3}}}^{3},3\,{z_{{3}}}^{3},\\
&& 0,-{z_{{3}}}^{2}z_{{28}}-z_{{7}}{z_{{3}}}^{2}-3\,z_{{8}}z_{{3}}z_{{4}}+3\,{z_{{4}}}^{2}+4\,z_{{5}}z_{{3}}-2\,{z_{{3}}}^{3}z_{{11}},-{z_{{3}}}^{3}z_{{13}}+9\,z_{{3}}z_{{4}}-2\,z_{{8}}{z_{{3}}}^{2},\\
&& -{z_{{3}}}^{3}z_{{14}}+3\,{z_{{3}}}^{2},6\,{z_{{3}}}^{2},0,0,10\,{z_{{3}}}^{3}z_{{4}},-z_{{7}}{z_{{3}}}^{3}-6\,z_{{8}}{z_{{3}}}^{2}z_{{4}}-2\,{z_{{3}}}^{3}z_{{28}}+10\,z_{{5}}{z_{{3}}}^{2}+15\,z_{{3}}{z_{{4}}}^{2},\\
&& 18\,{z_{{3}}}^{2}z_{{4}}-2\,z_{{8}}{z_{{3}}}^{3},3\,{z_{{3}}}^{3},45\,{z_{{3}}}^{2}{z_{{4}}}^{2}+20\,z_{{5}}{z_{{3}}}^{3},0,0,0,0] \\
X_{16}  &=& [0,0,0,0,0,0,0,0,0,0,0,0,0,0,0,0,0,0,0,0,0,0,0,0,0,0,0,0,0,0,0,0,0,0,1,0,0,\\
&& 0,0,0,0,0,0] \\
X_{17} &=&  [0,0,0,0,0,0,0,0,0,0,0,0,0,0,0,0,0,0,0,0,0,0,0,0,0,1,0,0,0,0,0,0,0,0,5\,z_{{3}},\\
&& -z_{{8}},3,0,15\,z_{{4}},0,0,0,0] \\
X_{18}  &=& [0,0,0,0,0,0,0,0,1,0,0,0,0,0,0,0,0,0,0,0,0,0,0,0,0,4\,z_{{3}},0,0,-z_{{8}},3,0,0,0,0,\\
&&
10\,{z_{{3}}}^{2},10\,z_{{4}}-4\,z_{{8}}z_{{3}},12\,z_{{3}},0,60\,z_{{3}}z_{{4}},0,0,0,0]
\end{eqnarray*}
\end{footnotesize}
\begin{footnotesize}
\begin{eqnarray*}
X_{19}  &=& [0,0,0,0,0,0,0,0,3\,z_{{3}},0,0,0,0,0,-z_{{8}},3,0,0,0,0,0,0,0,0,0,6\,{z_{{3}}}^{2},0,0,6\,z_{{4}}-3\,z_{{8}}z_{{3}},\\
&& 9\,z_{{3}},0,0,0,0,10\,{z_{{3}}}^{3},30\,z_{{3}}z_{{4}}-6\,z_{{8}}{z_{{3}}}^{2},18\,{z_{{3}}}^{2},0,90\,{z_{{3}}}^{2}z_{{4}},0,0,0,0] \\
X_{20} &=&  [0,0,0,0,0,0,0,0,3\,{z_{{3}}}^{2},0,0,0,0,0,3\,z_{{4}}-2\,z_{{8}}z_{{3}},6\,z_{{3}},0,0,0,0,0,0,0,0,0,4\,{z_{{3}}}^{3},0,0,\\
&& 12\,z_{{3}}z_{{4}}-3\,z_{{8}}{z_{{3}}}^{2},9\,{z_{{3}}}^{2},0,0,0,0,5\,{z_{{3}}}^{4},-4\,z_{{8}}{z_{{3}}}^{3}+30\,{z_{{3}}}^{2}z_{{4}},12\,{z_{{3}}}^{3},0,60\,{z_{{3}}}^{3}z_{{4}},0,0,0,0] \\
X_{21}  &=& [0,0,0,0,0,0,0,0,{z_{{3}}}^{3},0,0,0,0,0,-z_{{8}}{z_{{3}}}^{2}+3\,z_{{3}}z_{{4}},3\,{z_{{3}}}^{2},0,0,0,0,0,0,0,0,0,{z_{{3}}}^{4},0,0,\\
&& 6\,{z_{{3}}}^{2}z_{{4}}-z_{{8}}{z_{{3}}}^{3},3\,{z_{{3}}}^{3},0,0,0,0,{z_{{3}}}^{5},10\,{z_{{3}}}^{3}z_{{4}}-z_{{8}}{z_{{3}}}^{4},3\,{z_{{3}}}^{4},0,15\,{z_{{3}}}^{4}z_{{4}},0,0,0,0] \\
X_{22}  &=& [0,0,0,0,0,0,0,0,0,0,0,0,0,0,0,0,0,0,0,0,0,0,0,0,0,0,0,0,0,0,0,0,0,0,0,0,0,\\
&& 0,1,0,0,0,0] \\
X_{23} &=&  [0,0,0,0,0,0,0,0,0,0,0,0,0,0,0,0,0,0,0,0,0,0,0,0,0,0,0,0,0,0,0,0,0,0,0,1,0,\\
&& 0,6\,z_{{3}},0,0,0,0] \\
X_{24}  &=& [0,0,0,0,0,0,0,0,0,0,0,0,0,0,0,0,0,0,0,0,0,0,0,0,0,0,0,0,1,0,0,0,0,0,0,5\,z_{{3}},\\
&& 0,0,15\,{z_{{3}}}^{2},0,0,0,0] \\
X_{25}  &=& [0,0,0,0,0,0,0,0,0,0,0,0,0,0,1,0,0,0,0,0,0,0,0,0,0,0,0,0,4\,z_{{3}},0,0,0,0,0,0,\\
&& 10\,{z_{{3}}}^{2},0,0,20\,{z_{{3}}}^{3},0,0,0,0] \\
X_{26} &=&  [0,0,0,0,0,0,0,0,0,0,0,0,0,0,3\,z_{{3}},0,0,0,0,0,0,0,0,0,0,0,0,0,6\,{z_{{3}}}^{2},0,0,0,0,0,\\
&& 0,10\,{z_{{3}}}^{3},0,0,15\,{z_{{3}}}^{4},0,0,0,0] \\
X_{27}  &=& [0,0,0,0,0,0,0,0,0,0,0,0,0,0,3\,{z_{{3}}}^{2},0,0,0,0,0,0,0,0,0,0,0,0,0,4\,{z_{{3}}}^{3},0,0,0,0,0,\\
&& 0,5\,{z_{{3}}}^{4},0,0,6\,{z_{{3}}}^{5},0,0,0,0] \\
X_{28}  &=& [0,0,0,0,0,0,0,0,0,0,0,0,0,0,{z_{{3}}}^{3},0,0,0,0,0,0,0,0,0,0,0,0,0,{z_{{3}}}^{4},0,0,0,0,0,0,\\
&& {z_{{3}}}^{5},0,0,{z_{{3}}}^{6},0,0,0,0]
\end{eqnarray*}
\end{footnotesize}
\begin{footnotesize}
\begin{eqnarray*}
T_1 &=& [1,0,0,0,0,0,0,0,0,0,0,0,0,0,0,0,0,0,0,0,0,0,0,0,0,0,0,0,0,0,0,0,0,0,0,0,0,\\
&& 0,0,0,0,0,0] \\
T_2 &=& [0,0,-z_{{3}},-2\,z_{{4}},-3\,z_{{5}},-3\,z_{{6}},-2\,z_{{7}},-z_{{8}},-3\,z_{{9}},-2\,z_{{10}},-z_{{11}},-z_{{12}},0,z_{{14}},-3\,z_{{15}},\\
&& -2\,z_{{16}},-z_{{17}},-z_{{18}},0,z_{{20}},0,z_{{22}},2\,z_{{23}},3\,z_{{24}},-4\,z_{{25}},-4\,z_{{26}},-3\,z_{{27}},-2\,z_{{28}},-4\,z_{{29}},\\
&& -3\,z_{{30}},-2\,z_{{31}},-2\,z_{{32}},-z_{{33}},0,-5\,z_{{35}},-5\,z_{{36}},-4\,z_{{37}},-3\,z_{{38}},-6\,z_{{39}},z_{{40}},0,-z_{{42}},-2\,z_{{43}}] \\
T_3 &=& [0,0,0,-z_{{3}},-3\,z_{{4}},0,z_{{8}},-3,0,z_{{11}},0,2\,z_{{13}},z_{{14}},0,0,z_{{17}},0,2\,z_{{19}},z_{{20}},0,3\,z_{{22}},2\,z_{{23}},\\
&& z_{{24}},0,-6\,z_{{5}},-3\,z_{{6}},-2\,z_{{7}}+z_{{28}},-z_{{8}},-3\,z_{{9}},-2\,z_{{10}}+z_{{31}},-z_{{11}},2\,z_{{33}}-z_{{12}},z_{{34}},z_{{14}},\\
&& -10\,z_{{25}},-7\,z_{{26}},-5\,z_{{27}}+z_{{38}},-3\,z_{{28}},-15\,z_{{35}},0,0,-z_{{40}}z_{{3}},-z_{{42}}-2\,z_{{40}}z_{{4}}] \\
T_4 &=& [0,0,0,0,-z_{{3}},0,-1,0,0,0,0,0,0,0,0,0,0,0,0,0,0,0,0,0,-4\,z_{{4}},0,z_{{8}},-3,0,z_{{11}},\\
&& 0,2\,z_{{13}},z_{{14}},0,-10\,z_{{5}},-3\,z_{{6}},-2\,z_{{7}}+2\,z_{{28}},-z_{{8}},-20\,z_{{25}},0,0,0,-z_{{40}}z_{{3}}] \\
T_5 &=& [0,0,0,0,0,0,0,0,0,0,0,0,0,0,0,0,0,0,0,0,0,0,0,0,-z_{{3}},0,-1,0,0,0,0,0,0,0,\\
&& -5\,z_{{4}},0,z_{{8}},-3,-15\,z_{{5}},0,0,0,0] \\
T_6 &=& [0,0,0,0,0,0,0,0,0,0,0,0,0,0,0,0,0,0,0,0,0,0,0,0,0,0,0,0,0,0,0,0,0,0,-z_{{3}},\\
&& 0,-1,0,-6\,z_{{4}},0,0,0,0] \\
T_7 &=& [0,0,0,0,0,0,0,0,0,0,0,0,0,0,0,0,0,0,0,0,0,0,0,0,0,0,0,0,0,0,0,0,0,0,0,0,\\
&& 0,0,-z_{{3}},0,0,0,0]
\end{eqnarray*}
\end{footnotesize}
\section*{Appendix B: The
nonzero commutators for the Lie algebra $\mathcal{L}_{35}$ of the
differential operators of the homogeneous linear system of PDEs
(\ref{c3}) }
\begin{scriptsize}
\begin{equation}
\begin{tabular}{lllllll}
$[e_2, e_3] = e_2$,& $[e_2, e_5] = 2e_4$,& $[e_2, e_6] = e_5$,& $[e_2, e_8]= 3e_7$,& $[e_2, e_9] = 2e_8$,\\
$[e_2, e_{10}] = e_9$,& $[e_2, e_{12}] = 4e_{11}$,& $[e_2, e_{13}] = 3e_{12}$,&$[e_2, e_{14}] = 2e_{13}$,& $[e_2, e_{15}] = e_{14}$,\\
$[e_2, e_{17}] = 5e_{16}$,& $[e_2, e_{18}] = 4e_{17}$,& $[e_2, e_{19}] = 3e_{18}$,&$[e_2, e_{20}] = 2e_{19}$,& $[e_2, e_{21}] = e_{20}$,\\
$[e_2, e_{23}] = 6e_{22}$,& $[e_2, e_{24}] = 5e_{23}$,& $[e_2, e_{25}] = 4e_{24}$,&$[e_2, e_{26}] = 3e_{25}$,& $[e_2, e_{27}] = 2e_{26}$,\\
$[e_2, e_{28}] = e_{27}$,& $[e_2, e_{30}] = -e_2$,& $[e_2, e_{31}] = -e_4$,& $[e_2,e_{32}] = -e_7$,& $[e_2, e_{33}] = -e_{11}$,\\
$[e_2, e_{34}] = -e_{16}$,& $[e_2, e_{35}] = -e_{22}$,& $[e_3, e_4] = -e_4$,& $[e_3, e_6] = e_6$,& $[e_3, e_7] = -e_7$,\\
$[e_3, e_9] = e_9$,& $[e_3, e_{10}] = 2e_{10}$,& $[e_3, e_{11}] = -e_{11}$,& $[e_3,e_{13}] = e_{13}$,& $[e_3, e_{14}] = 2e_{14}$,\\
$[e_3, e_{15}] = 3e_{15}$,& $[e_3, e_{16}] = -e_{16}$,& $[e_3, e_{18}] = e_{18}$,& $[e_3,e_{19}] = 2e_{19}$,& $[e_3, e_{20}] = 3e_{20}$,\\
$[e_3, e_{21}] = 4e_{21}$,& $[e_3, e_{22}] = -e_{22}$,& $[e_3, e_{24}] = e_{24}$,& $[e_3,e_{25}] = 2e_{25}$,& $[e_3, e_{26}] = 3e_{26}$,\\
$[e_3, e_{27}] = 4e_{27}$,& $[e_3, e_{28}] = 5e_{28}$,& $[e_4, e_5] = 3e_7$,& $[e_4,e_6] = e_8$,& $[e_4, e_8] = 6e_{11}$,\\
$[e_4, e_9] = 3e_{12}$,& $[e_4, e_{10}] = e_{13}$,& $[e_4, e_{12}] = 10e_{16}$,& $[e_4,e_{13}] = 6e_{17}$,& $[e_4, e_{14}] = 3e_{18}$,\\
$[e_4, e_{15}] = e_{19}$,& $[e_4, e_{17}] = 15e_{22}$,& $[e_4, e_{18}] = 10e_{23}$,&$[e_4, e_{19}] = 6e_{24}$,& $[e_4, e_{20}] = 3e_{25}$,\\
$[e_4, e_{21}] = e_{26}$,& $[e_4, e_{30}] = -2e_4$,& $[e_4, e_{31}] = -3e_7$,& $[e_4,e_{32}] = -4e_{11}$,& $[e_4, e_{33}] = -5e_{16}$,\\
$[e_4, e_{34}] = -6e_{22}$,& $[e_5,e_6] = e_9$,& $[e_5, e_7] = -4e_{11}$,& $[e_5, e_9] = 2e_{13}$,& $[e_5, e_{10}] =2e_{14}$,\\
$[e_5, e_{11}] = -5e_{16}$,& $[e_5, e_{13}] = 3e_{18}$,& $[e_5, e_{14}] =4e_{19}$,& $[e_5, e_{15}] = 3e_{20}$,& $[e_5, e_{16}] = -6e_{22}$,\\
$[e_5, e_{18}] =4e_{24}$,& $[e_5, e_{19}] = 6e_{25}$,& $[e_5, e_{20}] = 6e_{26}$,& $[e_5, e_{21}] = 4e_{27}$,& $[e_5, e_{30}]= -e_5$,\\
$[e_5, e_{31}] = -e_8$,& $[e_5, e_{32}] = -e_{12}$,& $[e_5, e_{33}] = -e_{17}$,&$[e_5, e_{34}] = -e_{23}$,& $[e_6, e_7] = -e_{12}$,\\
$[e_6, e_8] = -e_{13}$,& $[e_6,e_{10}] = 2e_{15}$,& $[e_6, e_{11}] = -e_{17}$,& $[e_6, e_{12}] = -e_{18}$,& $[e_6, e_{14}]= 2e_{20}$,\\
$[e_6, e_{15}] = 5e_{21}$,& $[e_6, e_{16}] = -e_{23}$,& $[e_6, e_{17}] = -e_{24}$,& $[e_6,e_{19}] = 2e_{26}$,& $[e_6, e_{20}] = 5e_{27}$,\\
$[e_6, e_{21}] = 9e_{28}$,& $[e_7, e_8] = 10e_{16}$,& $[e_7, e_9] = 4e_{17}$,&$[e_7, e_{10}] = e_{18}$,& $[e_7, e_{12}] = 20e_{22}$,\\
$[e_7, e_{13}] = 10e_{23}$,&$[e_7, e_{14}] = 4e_{24}$,& $[e_7, e_{15}] = e_{25}$,& $[e_7, e_{30}] = -3e_7$,& $[e_7,e_{31}] = -6e_{11}$,\\
$[e_7, e_{32}] = -10e_{16}$,& $[e_7, e_{33}] = -15e_{22}$,& $[e_8,e_9] = 3e_{18}$,& $[e_8, e_{10}] = 2e_{19}$,& $[e_8, e_{11}] = -15e_{22}$,\\
$[e_8,e_{13}] = 6e_{24}$,&$[e_8, e_{14}] = 6e_{25}$,& $[e_8, e_{15}] = 3e_{26}$,& $[e_8, e_{30}] =-2e_8$,& $[e_8, e_{31}] = -3e_{12}$,\\
$[e_8, e_{32}] = -4e_{17}$,&$[e_8, e_{33}] = -5e_{23}$,& $[e_9, e_{10}] = 2e_{20}$,& $[e_9, e_{11}] = -5e_{23}$,&$[e_9, e_{12}] = -4e_{24}$,\\
$[e_9, e_{14}] = 4e_{26}$,&$[e_9, e_{15}] = 5e_{27}$,& $[e_9, e_{30}] = -e_9$,& $[e_9, e_{31}] = -e_{13}$,& $[e_9,e_{32}] = -e_{18}$,\\
$[e_9, e_{33}] = -e_{24}$,&$[e_{10}, e_{11}] = -e_{24}$,& $[e_{10}, e_{12}]= -2e_{25}$,& $[e_{10}, e_{13}] = -2e_{26}$,& $[e_{10}, e_{15}] = 5e_{28}$,\\
$[e_{11},e_{30}] = -4e_{11}$,&$[e_{11}, e_{31}] = -10e_{16}$,& $[e_{11}, e_{32}] = -20e_{22}$,& $[e_{12}, e_{30}] =-3e_{12}$,& $[e_{12}, e_{31}] = -6e_{17}$,\\
$[e_{12}, e_{32}] = -10e_{23}$,&$[e_{13}, e_{30}] = -2e_{13}$,& $[e_{13}, e_{31}] = -3e_{18}$,& $[e_{13}, e_{32}] =-4e_{24}$,& $[e_{14}, e_{30}] = -e_{14}$,\\
$[e_{14}, e_{31}] = -e_{19}$,&$[e_{14}, e_{32}] = -e_{25}$,& $[e_{16}, e_{30}] = -5e_{16}$,& $[e_{16}, e_{31}] =-15e_{22}$,& $[e_{17}, e_{30}] = -4e_{17}$,\\
$[e_{17}, e_{31}] = -10e_{23}$,&$[e_{18}, e_{30}] = -3e_{18}$,& $[e_{18}, e_{31}] = -6e_{24}$,& $[e_{19}, e_{30}] =-2e_{19}$,& $[e_{19}, e_{31}] = -3e_{25}$,\\
$[e_{20}, e_{30}] = -e_{20}$,&$[e_{20}, e_{31}] = -e_{26}$,& $[e_{22}, e_{30}] = -6e_{22}$,& $[e_{23}, e_{30}] = -5e_{23}$,&$[e_{24}, e_{30}] = -4e_{24}$,\\
$[e_{25}, e_{30}] = -3e_{25}$,&$[e_{26}, e_{30}] = -2e_{26}$,& $[e_{27}, e_{30}] = -e_{27}$,& $[e_{30}, e_{31}] = e_{31}$,&$[e_{30}, e_{32}] = 2e_{32}$,\\
$[e_{30}, e_{33}] = 3e_{33}$,&$[e_{30}, e_{34}] = 4e_{34}$,& $[e_{30}, e_{35}] = 5e_{35}$,& $[e_{31}, e_{32}] = 2e_{33}$,&$[e_{31}, e_{33}] = 5e_{34}$,\\
$[e_{31}, e_{34}] = 9e_{35}$,&$[e_{32}, e_{33}] = 5e_{35}.$&&&
\end{tabular}
\end{equation}
\end{scriptsize}

\end{document}